\documentclass[AMA,STIX1COL]{WileyNJD-v2}
\articletype{RESEARCH ARTICLE}%

\usepackage{float}
\usepackage{pgf,tikz}
\usepackage{pgfplots}
\usepackage{graphicx}   
\usepackage{longtable}   
\usepackage{epstopdf} 
\usepackage{pstcol,pst-fill,pst-grad}
\usepackage{subfigure}     
\usepackage{xcolor}
\newcommand{\ntext}[1]{\textcolor{blue}{#1}}

\newcommand{\LSS}{LSS}
\newcommand{\BLSS}{LSS\ }
\newcommand{\SLSS}{LSSs}
\newcommand{\BSLSS}{LSSs\ }
\newcommand{\MORPH}{\mathcal{S}}
\newcommand{\QNUM}{D}

\newcommand{\HYBINP}{\mathcal{U}}

\newcommand{\Rank}{\mathrm{rank}}
\newcommand{\AS}{\mathbf{S}}
\newcommand{\y}{\mathbf{y}}
\newcommand{\bu}{\mathbf{u}}
\newcommand{\WR}{span-reachable}

\newcommand{\SwitchSysLin}[1][{}]{ (n^{#1},Q,\{(A_{q}^{#1},B_{q}^{#1},C_{q}^{#1}) \mid q \in Q \},x_0^{#1})
}

\newcommand{\PARAMLinSwitch}[2]{ (n,Q,\{(A^{#2}_{q}({#1}),B^{#2}_{q}({#1}),C^{#2}_{q}({#1})) \mid q \in Q \},x_0^{#2}({#1})) }

\newcommand{\PARAMSET}{\Theta}
\DeclareMathOperator*{\MAPLSS}{\mathbf{\Pi^{\mbox{\tiny \textbf{LSS}}}}}
\newcommand{\MAPSARX}{\mathbf{\Pi^{\text{\tiny \textbf{SARX}}}}}

\renewcommand{\paragraph}[1]{\smallskip\noindent\textbf{#1.} }
\begin{document}

\title{Minimality and identifiability of discrete-time SARX systems}

\author[CL]{M. Petreczky*}

\author[CLb]{L. Bako}

\author[URIA]{S. Lecoeuche}

\author[CReSTIC]{K.M.D. Motchon}

\authormark{M. Petreczky \textsc{et al}}

\address[CL]{\orgdiv{CNRS}, \orgname{Centrale Lille}, \orgaddress{\state{UMR 9189--CRIStAL--Centre de Recherche en Informatique, Signal et Automatique de Lille}, \country{ F-59000 Lille, France}}}

\address[CLb]{\orgdiv{Ecole Centrale Lyon}, \orgname{Laboratoire Amp\`ere}, \orgaddress{\state{36 av. Guy de Collongue, 69130 Ecully}, \country{France}}}

\address[URIA]{\orgdiv{IMT Lille Douai}, \orgname{Univ. Lille}, \orgaddress{\state{Unit\'e de Recherche Informatique Automatique}, \country{ F-59000 Lille, France}}}

\address[CReSTIC]{\orgdiv{Université de Reims Champagne Ardenne}, \orgname{CReSTIC EA 3804}, \orgaddress{\state{51097 Reims}, \country{France}}}

%\address[2]{\orgdiv{Org Division}, \orgname{Org Name}, \orgaddress{\state{State name}, \country{Country name}}}

%\address[3]{\orgdiv{Org Division}, \orgname{Org Name}, \orgaddress{\state{State name}, \country{Country name}}}

\corres{M. Petreczky, CNRS, Centrale Lille, UMR 9189--CRIStAL--Centre de Recherche en Informatique, Signal et Automatique de Lille \email{mihaly.petreczky@centralelille.fr}}

%\fundingInfo{a FNADT subvention in the
%framework of the SUNRISE project.}

%\presentaddress{This is sample for present address text this is sample for present address text}

\abstract[Summary]{The paper studies the problem of minimality and identifiability for Switched Auto-Regressive eXogenous  (abbreviated by SARX) systems. We propose formal definitions  of the concepts of identifiability and minimality for SARX models. Based on these formalizations, we derive conditions for minimality and identifiability of SARX systems. In particular, we show that polynomially parametrized SARX systems are generically identifiable.}

\keywords{Switched ARX systems, minimality, identifiability, system identification, hybrid systems}

\jnlcitation{\cname{%
\author{M. Petreczky}, 
\author{L. Bako}, 
\author{S. Lecoeuche},and
\author{K.M.D. Motchon}} (\cyear{2019}), 
\ctitle{Minimality and identifiability of discrete-time SARX systems}, \cjournal{Int J Robust Nonlinear Control}, \cvol{2019;00:1--17}.}

\maketitle

%\footnotetext{\textbf{Abbreviations:} ANA, anti-nuclear antibodies; APC, antigen-presenting cells; IRF, interferon regulatory factor}

\section{Introduction}
System identification is the branch of control theory which is concerned with designing methods and algorithms for inferring parametrized mathematical models from input-output measurements. A fundamental criterion  characterizing the quality of a model parametrization is that of identifiability. This refers to the formal question of whether a given parametrized model can, in principle, be uniquely determined from  input-output data. More precisely, a parametrized model  structure is a map from a certain parameter space to a set of dynamic systems. Such a parametrized model structure is said to be (structurally) identifiable, if no two different parameter vectors yield two models whose input-output behavior is the same. The concept of identifiability has a number of implications for the design of informative experiments, the development of parameter estimation  algorithms, the  analysis of  identification methods and the significance of  estimated models. 

\paragraph{Contribution of the paper}
\color{blue}
The present paper deals with the problem of identifiability of switched ARX (\emph{abbreviated as SARX} ) systems. More precisely, we 
we introduce formal definitions identifiability for  SARX systems and show that a particular notion of minimality, called \emph{strong minimality},
is a sufficient condition for identifiability of SARX systems. We present conditions for checking strong minimality which are reminiscient of the well-know
minimality conditions for ARX systems. We also show that minimality that SARX parametrizations are generically minimal and generically identifiable. 
Note that minimality and identifiability are properties of the structure of the model  parametrization and not that of the data generated by the system. Based on our definitions, we derive checkable conditions guaranteeing these two properties. 
In addition to providing theoretical insights, the results of the paper allow us to check identifiability of SARX parametrizations, and to find identifiable parametrizations. 

It is worth noting that idenifiability and minimality  SARX systems cannot be reduced to the corresponding  properties of its ARX subsystems. 
It can be shown that a SARX system can be minimal \ntext{(resp. identifiable)}, even if none of the ARX subsystems is minimal (resp. identifiable). That is, the relationship between identifiability and minimality of SARX systems and their ARX subsystems is not straightforward.
\color{black} 
 
 %Our method of proof is roughly as follows: (1) convert SARX systems to a state-space form, using the regressors as a state-variable; (2) \ntext{Based on the results of Petreczky {\it et al}~\cite{MP:HSCC2010}}, we analyze minimality and identifiability of the resulting state-space representation. It should however be noted that the contributions of the paper are not trivial consequences of the work of Petreczky {\it et al}~\cite{MP:HSCC2010}. This  is due to the rich structure of  SARX systems which allows us to derive stronger results than for general linear switched systems. For example, the relationship between minimality and identifiability is much more direct for SARX systems than for linear switched systems. 

\paragraph{Motivation}
SARX systems are popular in the hybrid systems community, due to their simplicity and modelling power. In particular, most of hybrid systems identification algorithms were developed for SARX systems. Despite their popularity, identifiability and minimality of SARX systems are not yet completely understood. 

\color{blue}
Identifiability are essential for designing and analyzing algorithms for identification and adaptive control. Indeed, only identifiable parametrizations can be identified correctly by a parameter estimation algorithm. For this reason, identifiability is usually a necessary condition for correctness of parameter estimation algorithms. In turn, minimality is sufficient for achieving  identifiability of fully linearly parametrized model. 
% parameters which occur only in the non-minimal components of the system representation are not identifiable.

\color{black}
\paragraph{Related work}
Identification of hybrid systems is an active research topic \cite{Bako09-SYSID2,Bako08-IJC,VidalAutomatica,Verdult04, Juloski05-TAC, Roll04,Ferrari03, Nakada05,Paoletti2,Fliess08,IdentComparison, IdentSurvey, Garulli12-SYSID}. Many of the major contributions are formulated only for SARX systems \cite{Ozay12-TAC,Bako11-Automatica,Ohlsson13-Automatica,Lauer11-Automatica,VidalAutomatica,PaolettiTac}. The relationship between SARX systems and state-space representations was addressed by Paoletti {\it et al}~\cite{PaolettiTac} and Weiland {\it et al}~\cite{Weiland06}, and in this paper we use some of those results. 

To the best of our knowledge, the results of the paper are new. The paper of Vidal~\cite{VidalAutomatica} contains persistence of excitation conditions for SARX systems, which is related to identifiability. The main difference between the two  concepts is that the former  is a property of the data, while the latter is a property of the parameterization. Vidal~\cite{VidalAutomatica} also proposes a definition of minimality of SARX systems which implies our definition. However, the two definitions are not equivalent.

\ntext{Finally, we note that a preliminary version of the material presented in the current paper appeared in the proceedings of the 16th IFAC Symposium on System Identification~\cite{Petreczky2012-SYSID}. However, compared to the preliminary version, all detailed proofs of the SARX minimality and identifiability results are provided here.
They are mainly based on some technical results that give important properties of state space representation which arise from SARX systems. 
These technical results can be also used for analysising other structural properties of SARX systems.  
Moreover, we present an algorithm for finding an identifiable parametrization of SARX systems, the latter was not included in \cite{Petreczky2012-SYSID}. 
Furthermore, the organization of the paper has been improved in order to ease the reading. A motivating example is also provided.}

%The technical report \cite{SARXtech} represents an extended version of the current paper. 

\paragraph{Outline}
In Section~\ref{sect:sarx}, we define SARX systems and the corresponding system-theoretic concepts such as minimality and identifiability. In Section~\ref{sect:min}, we present sufficient conditions guaranteeing (strong) minimality of SARX systems. We also discuss the relationship between minimality of a SARX system and that of its subsystems. In Section~\ref{ident:sarx}, we discuss the relationship between minimality and identifiability. We show that minimality and identifiability are generic properties. Concluding remarks are provided in Section~\ref{sect:concl}.

%\cite{MP:HSCC2010,Petreczky2012-SYSID}

\paragraph{Notations}
Denote by $T=\mathbb{N}$ the time-axis of natural numbers. The notation described below is standard in formal language and automata theory~\cite{GecsPeak,AutoEilen}. Consider a set $X$ which will be called the \emph{alphabet}. Denote by $X^{*}$ the set of finite sequences of elements of $X$.  Finite sequences of elements of $X$ are  referred to as \emph{strings} or \emph{words} over $X$. We denote by $\epsilon$ the \emph{empty sequence (word)}. The length of a word $w$ is denoted by $|w|$, i.e. $|w|=k$ means that the length of $w$ is $k$; notice that $|\epsilon|=0$. We denote by $X^{+}$ the set of non-empty words, i.e. $X^{+}=X^{*}\setminus \{\epsilon\}$. We denote by $wv$ the concatenation of word $w \in X^{*}$ with $v \in X^{*}$ and recall that $\epsilon w=w\epsilon=w$. \ntext{Furthermore,  we denote by $I_{d}$ the $d \times d$ identity matrix and by $\textbf{O}_{d \times l}$ the $d \times l$ zero matrix.}

 %%%%%%%%%%%%%%%%%%%%%%%%%%%%%%%%%%%%%%%%%%%%%%%%%%%

\section{\ntext{Definition of SARX systems and problem formulation}}\label{sect:sarx}
\ntext{This section aims to introduce the concepts of minimality, strong minimality,  and identifiability for discrete-time SARX systems.  A formal description of this class of switched systems is recalled in Definition~\ref{SARX:def} 
}

\begin{Definition}[SARX systems]
 \label{SARX:def}
A SARX system $\AS$ of type $(n_y,n_u)$, where $0< n_u \le n_y$ are integers, is a collection
  \( \AS=\{h_q\}_{q \in Q} \), where \ntext{$Q$ is the finite set of discrete modes} and \ntext{for every $q\in Q$}, $h_q$ is a $p \times (n_y p+n_u m)$ matrices \ntext{with $p$ the output dimension and $m$ the input dimension of the system.} We will call a SARX system a SISO SARX system if $p=m=1$. \ntext{ The dimension of the SARX system $\AS$ is the number $pn_y+n_um$ and is denoted by $\dim \AS$.}
\end{Definition}

%%%%

Assigning semantics to SARX systems defined above requires that we first formalize the concept of input-output behaviour for SARX systems. \ntext{For this let's introduce the following notion of hybrid inputs of SARX systems.}
\color{blue}
\begin{Definition}[Hybrid inputs of SARX systems]
\label{notat:hybrid-inputs}
The hybrid inputs of SARX system $\AS$ in Definition~\ref{SARX:def} are the elements of $\HYBINP=Q \times \mathbb{R}^{m}$. For any $t\geq 0$, a sequence $w$ of the form
\begin{equation}
 \label{inp_seq}
 w=(q_0,u_0)\cdots (q_t,u_t) \in \HYBINP^{+}, 
 \end{equation}
where we recall that $\HYBINP^{+}$ denotes the set of non-empty finite sequences of elements of $\HYBINP$ describes the scenario, when discrete mode $q_i\in Q$ and continuous input $u_i\in\mathbb{R}^m$ are fed to $\AS$ at time $i$, for $i=0,\ldots,t$.
\end{Definition}
\begin{Notation}
\label{not:decomp}
In the sequel, \ntext{without loss of generality, we will set $Q=\{1,\ldots,\QNUM\}$} and we will use the following decomposition for the matrices $h_q$:
 \begin{eqnarray*}
    h_q&=&\begin{bmatrix} h_q^1 & h_q^2 & \cdots & h_q^{n_u+n_y} \end{bmatrix}, 
    %h_q^{1:n_y} &=& \begin{bmatrix} h_q^1, & \ldots, & h_q^{n_y}
   %               \end{bmatrix} \\
    %h_q^{n_y+1:n_y+n_u} &=& 
    %\begin{bmatrix} h_q^{n_y+1}, & \ldots, & h_{q}^{n_y+n_u}
    %\end{bmatrix}   
 \end{eqnarray*}
 where $h_q^{i} \in \mathbb{R}^{p \times p}$, $i=1,\ldots,n_y$,
 $h_{q}^{j} \in \mathbb{R}^{p \times m}$, $j=n_y+1,\ldots,n_u+n_y$.
\end{Notation} 
\color{blue} 
In order to introduce the formal definition of SARX minimality, the following concepts of input-output map realization as well as equivalence of SARX systems  are needed.
 %Now we are ready to define the semantics of SARX systems. 
\begin{Definition}[Input-output map realization and equivalence of SARX systems]
  %\begin{itemize}
  %\item
  %\item
  %$\hat{y}_i \in \mathbb{R}^{p}$, $i=1,\ldots, n_y$ are the
   %initial outputs.
  %\end{itemize}
  The SARX system $\AS$ is a realization of the input-output map $f:\HYBINP^{+} \rightarrow \mathbb{R}^{p}$, if for
  all $w \in \HYBINP^{+}$ of the form \eqref{inp_seq}, the outputs 
  \[ \y_i=f((q_0,\bu_0)\cdots (q_i,\bu_i)), i=0,\ldots,t \]
of $\AS$  satisfy the equation
  \begin{equation}
  \label{SARX:def1}
     \y_t= h_{q_t}\phi_t
  \end{equation}
  where we define the regressor 
  $\phi_t \in \mathbb{R}^{(n_yp+n_u m)}$ as
  \begin{equation} 
  \label{SARX:def2}
    \phi_t=\begin{bmatrix} \y_{t-1}^T & \y_{t-2}^T & \cdots & \y_{t-n_y}^T & \bu^T_{t-1} & \cdots \bu_{t-n_u}^T \end{bmatrix}^T,
  \end{equation}
 and for all $j < 0$, we set $\y_j=0$ and $\bu_j=0$. Two SARXs are called \emph{equivalent}, if they are realizations of the same input-output map.
 \end{Definition}
%\ntext{
For the input-output maps of interest $f:\HYBINP^{+} \rightarrow \mathbb{R}^{p}$, the value $f(w)$, with $w$ of the form \eqref{inp_seq}, describes the output of the system in Definition~\ref{SARX:def} at time $t$, generated as a response of the system to the hybrid input $w$. 
 Another important concept is that of minimality.
\begin{Definition}[Minimality of SARX systems]
 \label{sarx:def:min}
  A SARX system $\AS$ is minimal, if there exists no equivalent SARX of dimension less than
  $\dim \AS$.
 \end{Definition}
\color{blue}

Next, we define the concept of discrete-time \emph{linear switched system}~\cite{Sun:Book,Petreczky2012-Automatica} (abbreviated by \LSS) associated with a SARX.
This will be used to define the concept of strong  minimality, which will play a central role in the rest of the paper. 
We will use the notation and terminology from \cite{Petreczky2012-Automatica} for linear switched systems in state-space form, which we recal below.
% in Appendix~\ref{AppdxLSS} (cf. Definition~\ref{switch_sys:real:def1}). 
%Note also that input-output maps realization of SARX systems can be also characterized by input-output maps realization of their
 %associated 
%
% see Lemma~\ref{sarx2lss:lemma0} below. 
%
%These \LSS\, are given in Definition~\ref{sarx2lss}. The input-output map realization for  \LSS\,. See for instance the work of Petreczky {\it et al}~\cite{MP:HSCC2010} for more details on this concept.
\begin{Definition}%{(Linear Switched Systems)}
\label{lin_switch:def}
A \emph{linear switched system} (abbreviated by \LSS) is a discrete-time system $\Sigma$ represented by
\begin{equation}
\label{lin_switch0}
%\Sigma: \left \{ 
  \begin{array}{lcl}
   x_{t+1}&=& A_{q_t}x_t+B_{q_t}u_t,\, x_0 \mbox{ fixed}\\
   y_t&=&C_{q_t}x_t.
\end{array}
%\right.
\end{equation}
where, 
 $x_t \in \mathbb{R}^{n}$ is the continuous state at time $t \in T$,
%\item
 $u_t \in \mathbb{R}^{m}$ 
 is the continuous input at time $t \in T$, 
 $y_t \in \mathbb{R}^{p}$ is the continuous output at time $t \in T$,
%\item 
$q_t \in Q$ is the discrete mode (state) at time $t$,  $Q$ is the finite
 set of discrete modes, and
 $x_0 \in \mathbb{R}^n$ is the initial state of $\Sigma$.
For each discrete mode $q \in Q$,
the corresponding matrices are of the form
$A_{q} \in \mathbb{R}^{n \times n}$, 
$B_{q} \in \mathbb{R}^{n \times m}$ and $C_{q} \in \mathbb{R}^{p \times n}$. 
\end{Definition}
%%%%\vspace{-20pt}
\begin{Notation}
We will use $\SwitchSysLin$ as a short-hand notation for \SLSS\  of the form (\ref{lin_switch0}).
\end{Notation}
\begin{Definition}[Associated \LSS\,of SARX]%[\cite{Weiland06}]
 \label{sarx2lss}
 Let $\AS=\{h_q\}_{q \in Q}$ be a SARX system of type $(n_y,n_u)$. 
 The \LSS\ $\Sigma_{\AS}=\SwitchSysLin$ associated with the SARX $\AS$ is given by: 
%\begin{equation}
%\label{lin_switch01}
%\Sigma: \left \{ 
%  \begin{array}{lcl}
%   x_{t+1}&=& A_{q_t}x_t+B_{q_t}u_t,\, x_0 \mbox{ fixed},\\
%   y_t&=&C_{q_t}x_t,
%\end{array}
%\right.
%\end{equation}
%where, 
% $x_t \in \mathbb{R}^{n}$ with $n=pn_y+mn_u$ is the continuous state at time $t \in T$,
%\item
% $u_t \in \mathbb{R}^{m}$ 
% is the continuous input at time $t \in T$, 
% $y_t \in \mathbb{R}^{p}$ is the continuous output at time $t \in T$,
%\item 
%$q_t \in Q$ is the discrete mode (state) at time $t$,
% $x_0 \in \mathbb{R}^n$ is the initial state of $\Sigma_{\AS}$ and 
 \begin{equation}
 \label{sarx2lss:eq1}
\begin{array}{lll}
A_q = 
\begin{bmatrix}
A_q^{y} & A_{q}^{u}
\end{bmatrix}
,\quad\quad
C_q=h_q,
\quad\quad
B_q=
\begin{bmatrix} 
\mathbf{O}_{pn_y \times m} \\
I_{m} \\ 
\mathbf{O}_{m(n_u-1) \times m}
\end{bmatrix}, \quad x_0=0 
\end{array}
\end{equation} 
with 
\begin{equation}
 \label{sarx2lss:eq1bis}
A_{q}^{y}=
\begin{bmatrix}
\begin{bmatrix} h_{q}^1 & h_{q}^{2} & \ldots & h_{q}^{n_y-1} \end{bmatrix} & h_q^{n_y} \\
I_{(n_y-1)p} & \mathbf{O}_{(n_y-1)p \times p}\\
\mathbf{O}_{m \times p(n_y-1)} & \mathbf{O}_{m \times p} \\ 
\mathbf{O}_{(n_u-1)m \times p(n_y-1)} & \mathbf{O}_{(n_u-1)m \times p}
\end{bmatrix},
\quad\quad 
A_{q}^{u}
= 
\begin{bmatrix}
\begin{bmatrix} h_{q}^{n_y+1} & h_{q}^{n_y+2} & \cdots & h_{q}^{n_y+n_u-1} \end{bmatrix} & h_q^{n_u+n_y} \\
\mathbf{O}_{(n_y-1)p \times (n_u-1)m} & \mathbf{O}_{(n_y-1)p \times m} \\
\mathbf{O}_{m \times m(n_u-1)} & \mathbf{O}_{m \times m} \\
I_{(n_u-1)m} & \mathbf{O}_{(n_u-1)m \times m}
\end{bmatrix},
\end{equation}
and $h_q$ decomposed as in Notation \ref{not:decomp}. %For clarity of the presentation the associated \LSS\,~\eqref{lin_switch0}--\eqref{sarx2lss:eq1bis} with the SARX system $\AS=\{h_q\}_{q \in Q}$ is also denoted by \( \Sigma_{\AS}:=\SwitchSysLin \). %The dimension of $ \Sigma_{\AS}$ is the dimension $n$ of its continuous state-space.
 \end{Definition}
Similarly to the case of SARX systems, we can define the concept of a \LSS\ being a realization of an input-output map $f:\HYBINP^{+} \rightarrow \mathbb{R}^{p}$. 
% see Appendix~\ref{AppdxLSS}. 
Informally, a \LSS\ is a realization of $f$, if the for any sequence of discrete modes and inputs $w \in \HYBINP$,
the output response of the \LSS\ from its initial state to this sequence of  inputs and discrete equals $f(w)$. 
Formally, consider a state $x_{init} \in \mathbb{R}^{n}$. For any input sequence $w \in \HYBINP^{*}$, 
     let $x_{\Sigma}(x_{init},w)$ be the state of $\Sigma$
     reached from $x_{init}$ under input $w$, % at time $t=|w|$. 
     i.e.
     $x_{\Sigma}(x_{init},w)$ is defined recursively as follows;
     $x_{\Sigma}(x_{init},\epsilon)=x_{init}$,  and if $w=v(q,u)$
     for some $(q,u) \in \HYBINP$, $v \in \HYBINP^{*}$, then
     \( x_{\Sigma}(x_{init},w)=A_{q}x_{\Sigma}(x_{init},v)+B_{q}u. \)
     If $w \in \HYBINP^{+}$, then denote
      by $y_{\Sigma}(x_{init},w)$ the 
      output response of $\Sigma$ to $w$, from the state 
      $x_{init}$, i.e.
      %state $x$ if input $w$ is fed to $\Sigma$. That is,
      if $w=v(q,u)$ for some $(q,u) \in \HYBINP$, $v \in \HYBINP^{*}$,
      then 
      \( y_{\Sigma}(x_{init},w)=C_{q}x_{\Sigma}(x_{init},v). \)
      The function  
       \begin{displaymath}
       y_{\Sigma}: \HYBINP^{+} \rightarrow \mathbb{R}^{p},\quad w\longmapsto y_{\Sigma}(w)=y(x_0,w),
\end{displaymath} 
is called the  \emph{input-output map} of $\Sigma$.   
%The input-output map of $\Sigma$ maps each sequence $w \in \HYBINP^{+}$ to the output generated by $\Sigma$ if $w$ is fed to $\Sigma$ starting from a given initial state $x_0$. The definition above implies that the input-output behavior of a \BLSS can be described as a map $f:\HYBINP^{+} \rightarrow \mathbb{R}^{p}$.
 %%\vspace{-15pt}
%   \begin{Definition}[\SLSS\,realization of input-output maps]
%   \label{switch_sys:real:def1}
    An input-output map $f: \HYBINP^{+} \rightarrow \mathbb{R}^{p}$
      is said to be \emph{realized} by a \LSS\ $\Sigma$  of the form (\ref{lin_switch0}) if $f$ equals the input-output map  $y_{\Sigma}$ of $\Sigma$.
   %\[ \forall w \in \HYBINP^{+},\quad  f(w)=y_{\Sigma}(x_0,w). \]
    In this case $\Sigma$ is said to be a realization of $f$. We say that two \LSS{s} are \emph{equivalent}, if their input-output maps are equal. 
\begin{Lemma}
\label{sarx2lss:lemma0}
 The SARX system $\AS$ is a realization of the input-output map $f$ if and only if the associated \LSS\ $\Sigma_{\AS}$ is a realization of $f$.
\end{Lemma}
\begin{proof}
See the paper of Weiland {\it et al}~\cite{Weiland06}.
\end{proof}
The following corollary of Lemma~\ref{sarx2lss:lemma0} allows us to relate the problem of minimality of SARX to that of \SLSS.  The latter has already been investigated by Petreczky {\it et al}~\cite{Petreczky2012-Automatica}.
Formally,  the dimension of  a LSS $\Sigma$ of the form \eqref{lin_switch0}, denoted by $\dim \Sigma$, is defined as the dimension $n$ of its continuous state-space.  Let $f:\HYBINP^{+} \rightarrow \mathbb{R}^{p}$ be an
 input-output map and let $\Sigma$ be a \LSS\  which is a realization of $f$.
 Then $\Sigma$ is \emph{a minimal realization of $f$}, if for any \LSS\  realization
 $\hat{\Sigma}$ of $f$, $\dim \Sigma \le \dim \hat{\Sigma}$.
That is, a \LSS\  realization is a minimal realization of $f$
if it has the smallest dimension of state-space among all 
the \BLSS which are realizations of $f$. We will say that a \LSS\ $\Sigma$ is a \emph{minimal}, if it is a minimal realization of its own input-output map $y_{\Sigma}$. 
% A formal definition \SLSS\, minimality is also recalled in Appendix~\ref{AppdxLSS} (see Definition).
 \begin{Corollary}
 \label{sarx2lss:col1}
  If the associated \LSS\ $\Sigma_{\AS}$ is minimal, then $\AS$ is minimal. 
 \end{Corollary}
%%%%%%%%%%%%%%%%%%%%%%%%
 \begin{proof}[Proof of Corollary \ref{sarx2lss:col1}]
  Assume that $\AS$ is not minimal. Then there exists an
  equivalent $\AS_{m}$ of type $(n_y^{'},n_u^{'})$ such that
  $n_{y}^{'}+n_{u}^{'} < n_y+n_u$. But this implies that
  $\dim \Sigma_{\AS_{m}}=n_{y}^{'}+n_{u}^{'} <  n_y+n_u=\dim \Sigma_{\AS}$, which contradicts to the minimality of
  $\Sigma_{\AS}$.
 \end{proof}
\begin{remark}
It can be noticed that none of the linear subsystems of
$\Sigma_{\AS}$ in \eqref{sarx2lss:eq1}--\eqref{sarx2lss:eq1bis} is minimal. Indeed,  for each $q \in Q$, $A_q$ contains a zero
  row, hence $\Rank A_{q} < n_u+n_y$. This means that
  $\lambda=0$ is an eigenvalue of $A_q$. By the PBH criterion, $(C_q,A_q)$ is an observable pair
  if and only if the matrix
  $\big[\begin{matrix} C_q^T, & \lambda I-A_q^T \end{matrix}\big]$
  has rank $n_y+n_u$ for all the eigenvalues of $A_q$.
  We will show that for $\lambda=0$ this matrix cannot
  be of full row rank. To see this, for $\lambda=0$ the matrix
  becomes $\big[\begin{matrix} C_q^T, -A_{q}^T \end{matrix}\big]$.
  But $C_q$ equals the first row of $A_q$ multiplied
  by $-1$. Hence, $[\begin{matrix} C_q^T, -A_{q}^T \end{matrix}]$ will have the same rank as $A_q$ and that is smaller
  than $n_u+n_y$. Thus, the linear subsystems  $\Sigma_{\AS}$ in \eqref{sarx2lss:eq1}--\eqref{sarx2lss:eq1bis} are not observable and consequently, they are not minimal. However, as we shall see later, they are generically minimal.
\end{remark}
 The result of Corollary~\ref{sarx2lss:col1} prompts us to propose the following definition.
\color{black}
 \begin{Definition}[Strong minimality of SARX systems]
  A SARX system $\AS$ is called \emph{strongly minimal}, if the corresponding \LSS\ $\Sigma_{\AS}$ is minimal.
 \end{Definition}
 %Note that in the definition of a SARX involves specifying the 
 %initial state which is the collection of the first $n_y$ outputs.

\begin{remark}
\label{example:min1}
By Corollary \ref{sarx2lss:col1}, strong minimality implies minimality. However, Minimality does not imply strong minimality. Indeed, consider the SARX system $\AS$ with discrete modes $Q=\{1,2\}$ such that the ARX subsystem associated with mode $1$ is $\y_t=-\y_{t-2}+\bu_{t-1}$ and the ARX subsystem associated with mode $2$ is $\y_t=-2\y_{t-2}+2\bu_{t-1}$. The two ARX subsystems are distinct, each of them is minimal, yet the associated \LSS\ $\Sigma_{\AS}$ is not minimal (in fact, it is not observable). The latter can be checked using the minimality conditions of Petreczky {\it et al}~\cite{MP:HSCC2010}.
\end{remark} 
 
\begin{remark}
\label{example:min2}
Minimality of the ARX subsystems is not necessary for strong minimality (and hence minimality) of the whole system. To see this, consider again the SARX system $\AS$ with two discrete modes $Q=\{1,2\}$ such that the ARX subsystem in mode $1$ is of the form $\y_{t}=8\y_{t-1}-15\y_{t-2}+\bu_{t-1}-3\bu_{t-2}$, and the ARX subsystem in mode $2$ is of the form $\y_{t}=\y_{t-1}+2\y_{t-2}+\bu_{t-1}+\bu_{t-2}$. The transfer function of the ARX in the first mode is $\frac{z-3}{z^2-8z+15}=\frac{1}{z-5}$ and the transfer function of the second ARX is $\frac{z+1}{z^2-z-2}=\frac{1}{z-2}$, hence neither of them is minimal. Yet, by using the conditions of Petreczky {\it et al}~\cite{MP:HSCC2010}, it can be easily shown that the \LSS\ $\Sigma_{\AS}$ is minimal. Since strong minimality implies minimality of SARX systems, we get that $\AS$ is minimal. 
\end{remark}

In order to be able to speak of identifiability, we need the notion of parametrization of SARX systems.  
 \begin{Notation}
  Denote by $SARX(n_y,n_u,m,p,Q)$ the set of all SARX systems of
  type $(n_y,n_u)$ with input space $\mathbb{R}^{m}$,
  output space $\mathbb{R}^{p}$, and set of discrete modes $Q$.
 \end{Notation}
%%%%\vspace{-20pt}
 \begin{Definition}{(Parametrization of SARX systems)}
 \label{sarx:param:def}
 Assume that $\PARAMSET \subseteq \mathbb{R}^{d}$ is the set of parameters.
 A \emph{SARX parametrization} is a map
%belonging to $SARX(n_y,n_u,m,p,Q)$} is a map
 \begin{equation}
 \label{sarx:param:def:eq1}
    \MAPSARX:\PARAMSET \rightarrow SARX(n_y,n_u,m,p,Q) 
  \end{equation}
  %For each $\theta \in \PARAMSET$,
  %we denote $\MAPSARX(\theta)$ by
  %\[ \AS(\theta)=(\{h_q(\theta) \}_{q \in Q})  \]
 \end{Definition}
\color{blue}
\begin{Example}\label{Exple_Pbme}
Consider the discrete-time model of the intake manifold of a spark ignition engine as described in \cite{Toth12-CST} and \cite{Kwiatkowski-PhD2008}, 
$$y_t=f(u_t,p_t,v_t),$$ where 
the output $y_t$ is the normalized air charge, the input $u_t$ is the opening of the throttle valve;  $p_t$ and $v_t$ refer respectively to the pressure inside the intake manifold and the speed of the engine. Here, $p$ and $v$ are viewed as external signals which take values in some bounded intervals.  
We refer to \cite{Toth12-CST} and \cite{Kwiatkowski-PhD2008}  for more details. 
Inspired by the Linear Parameter Varying (LPV) parameterization of the intake model proposed in \cite{Toth12-CST}, one  can consider approximating the intake manifold  with a SISO SARX system of type $(2,2)$ by viewing $p$ and $v$ as piecewise constant signals, each of which is allowed, for simplicity, to take only two possible values, $(p_1,p_2)=(0.3,0.7)$,  and $(v_1,v_2)=(2,5)$ respectively. Indeed, since $p$ and $v$ are external excitation signals, one can select them in such a way that they are exactly piecewise constant with values prescribed as above. In that case, the SARX model considered here coincides with the LPV one. 

 Let us form a vector  $\bar{p}_{(i,j)}=\begin{bmatrix}10^{-3} & p_i & p_i^2 & v_j & v_j^2\end{bmatrix}^\top $ and introduce the bijective map $\sigma:\left\{1,2\right\}\times \left\{1,2\right\} \rightarrow Q=\left\{1,2,3,4\right\}$, $(i,j)\mapsto q=\sigma(i,j)$. 
where $p_1=0.3$, $p_2=0.7$, $v_1=2$, $v_2=5$. 
Then the LPV parameterization in \cite{Toth12-CST} reduces to a SARX one $\MAPSARX: \Theta=\Re^{20} \rightarrow SARX(2,2,1,1,Q) $ with $\MAPSARX(\theta)=\left\{h_q(\theta)\right\}_{q\in Q}$,  
$ h_q(\theta)=M(\theta)\bar{p}_{(i,j)}$ for $q=\sigma(i,j)$ and  
$$M(\theta)= \begin{bmatrix}
	\theta_{1} &\theta_{5}& \cdots & \theta_{17}\\
	\theta_{2} & \theta_{6} &\ldots & \theta_{18}\\
	\theta_{3}&  \theta_{7} & \cdots & \theta_{19}\\
	\theta_{4} & \theta_{8} & \ldots & \theta_{20}
\end{bmatrix}\in \mathbb{R}^{4\times 5}$$
Here, $\bar{p}_{(i,j)}$ is a vector of the form 
$\bar{p}_{(i,j)}=\begin{bmatrix}10^{-3} & p_i & p_i^2 & v_j & v_j^2\end{bmatrix}^\top $
where $p_1=0.3$, $p_2=0.7$, $v_1=2$, $v_2=5$. 
$p_1\in \left\{311.75;737.25\right\}$ and $p_2\in \left\{2132.5;4877.5\right\}$
he input-output map of such a system is hence defined by 
$$y_t=h_{q_t}^\top \phi_t $$
with $\phi_t=\begin{bmatrix}y_{t-1}& y_{t-2} & u_{t-1} & u_{t-2}\end{bmatrix}^\top$. 
This example illustrates how the behavior of some physical systems can be modelled by SARX systems.
% Indeed, this form of model  is largely used for switched system identification, see e.g., \cite{Paoletti07,Garulli12}.
\end{Example}
\color{black}
 \begin{Definition}[Identifiability of SARX parametrizations] The parametrization $\MAPSARX$ is called \emph{identifiable}, if for $\theta_1 \ne \theta_2 \in \PARAMSET$, the corresponding SARX $\MAPSARX(\theta_1)$ and $\MAPSARX(\theta_2)$ are not equivalent.
 \end{Definition}
\noindent \ntext{The intuition behind the above definition is that if a parametrization is not identifiable, then there might exist different parameter values which yield the same observed behaviour and hence they cannot be distinguished from each other by input-output experiments.
}
\color{blue}
Hence, the problem of identifying the parameters of a SARX models from a non identifiable parametrization is ill-posed. It would be tempting to try to reduce identifiability of SARX parametrizations to that of the parametrization of the corresponding ARX subsystems. This would then allow us to use existing theory on identifiability of ARX parametrizations. 
\textbf{Unfortunately, identifiability of a SARX parametrization does not imply the 
 identifiability of the corresponding parametrization  of
 ARX subsystems.} The example below demonstrates this point.
 \begin{Example}
 \label{ident:example}
  Consider the SARX parametrization $\MAPSARX$ with $\PARAMSET=\mathbb{R}^{2}$, 
  and consider the parametrization 
  $\MAPSARX((\theta_1,\theta_2))=\{h_{q}(\theta_1,\theta_2)\}_{q \in Q}$, where 
  %$h_{1}=\big[\begin{matrix} (\theta_1+\theta_2) & -\theta_1\theta_2, & 1  & -\theta_2 \end{matrix}\big]$ and $h_{2}=\big[\begin{matrix} (2+\theta_2) & -2\theta_2, & 1  & -\theta_2 \end{matrix}\big]$.
  \begin{displaymath}
  h_{1}=\big[\begin{matrix} (\theta_1+\theta_2) & -\theta_1\theta_2 & 1  & -\theta_2 \end{matrix}\big],
  \end{displaymath}
  and
  \begin{displaymath}
   h_{2}=\big[\begin{matrix} (2+\theta_2) & -2\theta_2 & 1  & -\theta_2 \end{matrix}\big].
  \end{displaymath}
%   \[ 
%     \begin{split}
%      & h_{1}=\big[\begin{matrix} (\theta_1+\theta_2) & -\theta_1\theta_2, & 1  & -\theta_2 \end{matrix}\big] \\
%    &  h_{2}=\big[\begin{matrix} (2+\theta_2) & -2\theta_2, & 1  & -\theta_2 \end{matrix}\big] 
%    \end{split}
%   \]
   Define the set 
  $G=\{(\theta_1,\theta_2) \mid \theta_1 \ne 2\}$.
  Consider the restriction $\MAPSARX|_{G}$ of 
  $\MAPSARX$ to $G$. Using Theorem \ref{sarx2lss:lemma3} 
  one can check
   that for
  any $(\theta_1,\theta_2) \in G$, the SARX system
  $\MAPSARX((\theta_1,\theta_2))$ is strongly minimal.
  Hence, the parametrization $\MAPSARX|_{G}$ is identifiable
  by Theorem \ref{sarx:main2.1}. Identifiability of $\MAPSARX|_{G}$ can
  also be checked by considering the
  switching sequence $112$ and input $\bu_0=1$, $\bu_t=0$, $t > 0$ and noticing that 
  then $\y_0$, $\y_1=1$, $\y_2=\theta_1$, 
  $\y_3=2\theta_1+\theta_2\theta_1-2\theta_2$ from which
  $\theta_2=\frac{(\y_3-2\theta_1)}{\theta_1-2}$.
  Hence, $\theta_1$ and $\theta_2$ can be determined from
  the outputs $\y_2$ and $\y_3$. \\
  Note however, that for any $(\theta_1,\theta_2)$, the
   ARX subsystems of $\MAPSARX(\theta_1,\theta_2)$ are not
   identifiable, since their dynamics does not depend on 
   $\theta_2$.
 \end{Example}
This implies that identifiability of SARX parametrizations has to be investigated separately.  
Recall now that identifiability of ARX parametrizations is closely related to their minimality. Hence, we start by investigating mnimality of SARX models.
\color{black}

\section{Main results}
 \textcolor{blue}{In this section we present the main results of the paper. First, in Section \ref{sect:min} we discuss minimality of SARX systems. In Section \ref{ident:sarx} we use the
 results of Section \ref{sect:min} to characterize identifiability of SARX systems.}

\subsection{Minimality conditions for switched ARX systems}
\label{sect:min}
In this section we will analyze minimality of SARX systems.\ntext{ We start by Lemmas~\ref{element:lemma1},~\ref{element:lemma3} and~\ref{element:lemma4} that provide some simple and crucial properties of minimal SARX systems}
 \begin{Lemma}
 \label{element:lemma1}
  If the SISO SARX system $\AS$ is minimal, then there must exist $q\in Q$ such that
  $h_{q}^{n_u+n_y} \ne 0$.
 \end{Lemma}
%%%%%%%%%%%%%%%%%%%%%%%%%
 \begin{proof}%[Proof of Lemma \ref{element:lemma1}]
  Assume the contrary, i.e. that $h_{q}^{n_u+n_y}=0$ for
  all $q \in Q$. Define the vectors $\hat{h}_q \in \mathbb{R}^{n_y+n_u-1}$ by 
  $\hat{h}_q=\big[\begin{matrix} h_q^1 & \ldots & h_q^{n_u+n_y-1}\end{matrix}\big]$, $q \in Q$. Define the regressors $\hat{\phi}_t$ as 
  $$\hat{\phi}_t=\big[\begin{matrix} \y_{t-1}^T & \ldots, \y_{t-n_y}^T & \bu_{t-1}^T & \ldots &\bu_{t-(n_u-1)}^T \end{matrix}\big]^T,$$
  where we used the convention that
  $\y_{j}=0$ and $\bu_{j}=0$ for $j < 0$. It then follows that
  $\y_t=h_{q}\phi_t=\hat{h}_q\hat{\phi}_t$
  for all $t \in T$. Hence, $\hat{\AS}=(\{\hat{h}_q\}_{q \in Q})$  realizes the same input-output map as $\AS$.
  But the dimension of $\hat{\AS}$ is smaller than that of
  $\AS$, which contradicts the minimality of $\AS$.
 \end{proof}

\begin{Lemma}
 \label{element:lemma3}
A SISO ARX system is minimal according to Definition \ref{sarx:def:min}  if and only if the numerator and denominator of its transfer function are co-prime.
 \end{Lemma}
\begin{proof}%[Proof of Lemma \ref{element:lemma3}]
The proof follows from the classical linear theory, by observing that two ARX systems realize the same input-output map if and only if they have the same transfer function (modulo zero/pole cancellation).

Consider an ARX system $\y_t=h_{q_t}\phi_t$ and assume that it is minimal. If its transfer function admits a zero-pole cancellation, then the degrees of the numerator and denominator of the transfer function decrease by one. The latter means that the transfer function can be realized by an ARX of type $(n_y-1,n_u-1)$. The dimension of the latter is $n_y+n_u-2$ and hence smaller than that of the original system, which was supposed to be minimal. Moreover, this new ARX system will realize the same input-output map as the original one.

Conversely, consider an ARX system $\AS$ whose transfer function does not allow zero/pole cancellation.  Let $f$ be the input-output map of $\AS$ and assume that the ARX system $\hat{\AS}$ is a minimal realization of $f$. Then the transfer function $H_{\hat{\AS}}(z)$ cannot allow a zero/pole cancellation and it must be equal to the transfer function $H_{\AS}(z)$ of $\AS$. Since neither $H_{\AS}(z)$ nor $H_{\hat{\AS}}(z)$ allow zero/pole cancellation, their equality implies the equality of the numerators and denominators respectively, viewed as polynomials. In particular, the corresponding coefficients are the same and hence the parameters of the two ARX systems are the same too. In particular, the dimensions of the two systems will be the same, and hence $\AS$ is then a minimal realization of its input-output map.
 \end{proof}
\ntext{
\begin{remark}
Recall that in the classical literature, a SISO ARX is said to be minimal if and only if the numerator and the denominator of its transfer function are co-prime polynomials. Consequently, Lemma~\ref{element:lemma3} shows that our definition of minimality is consistent with the traditional one. 
\end{remark}
}

\begin{Lemma}
 \label{element:lemma4}
If at least one of the ARX subsystems of a SISO SARX system is minimal, then the system is minimal.
\end{Lemma}
%%%%%%%%%%%%%%%%
\begin{proof}%[Proof of Lemma \ref{element:lemma4}]
  Consider $\AS=\{h_q\}_{q \in Q}$ and assume that for some $q_m \in Q$, the ARX $\y_t=h_{q_m}\phi_t$ is minimal. Assume that $\AS$ is not minimal and hence there exists a SARX $\AS_m=(\{\hat{h}_q\}_{q \in Q})$ such that $\dim \AS_m \le \dim \AS$ and $\AS_m$ realizes the same input-output map as $\AS$. It then follows that the dimension of the ARX $\y_t=h_{q_m}\phi_t$ is 
larger than that of $\y_t=\hat{h}_{q_m}\phi_t$. It also follows that both $\y_t=\hat{h}_{q_m}\phi_t$ and $\y_t=h_{q_m}\phi_t$ realizes the same linear input-output map\footnote{Namely, the map which maps input $\bu_0,\ldots,\bu_{t}$ to the output $\y_t=f((q_m,\bu_0)\cdots (q_m,\bu_{t}))$, where $f$ is the input-output map of $\AS$}. This contradicts the minimality of $\y_t=h_{q_m}\phi_t$.
 \end{proof}

The definition \ref{sarx:def:min} of minimality for SARX systems might seem ambiguous because it does not exclude explicitly the possibility of having two minimal SARX realizations of types $(n_y,n_u)$ and $(\hat{n}_y,\hat{n}_u)$ respectively for the same input-output map with $(n_y,n_u) \ne (\hat{n}_y,\hat{n}_u)$. According to the lemma below, this is impossible at least in the SISO case.
\begin{Lemma}
 \label{element:lemma5}
Assume that $\AS_1$ and $\AS_2$ are two minimal and equivalent SISO SARX systems such that $\AS_1$ is of type $(n_y,n_u)$ and $\AS_2$ is of
  type $(\hat{n}_y,\hat{n}_u)$. Then $(n_y,n_u)=(\hat{n}_y,\hat{n}_u)$.
 \end{Lemma}
%%%%%%%%%%%%%%%%%%%%%%%%%%%%%%%%%%
 \begin{proof}%[Proof of Lemma \ref{element:lemma5}]
Pick any discrete state $q$ and consider the transfer functions $H_i(z)$, $i=1,2$, of the ARX system in mode $q$ associated with the SARX $\AS_i$, $i=1,2$.  Since $\AS_1$ and $\AS_2$ are equivalent, they produce the same response to any input if the discrete mode is kept to be $q$. Hence, the ARX systems corresponding to the mode $q$ are also equivalent, i.e. $H_1(z)$ and $H_2(z)$ describe the same input-output behavior. This means that the transfer functions $H_1(z)$ and $H_2(z)$ are equal as rational expressions, after possibly performing zero/pole cancellation. The degrees of the numerators of $H_1(z)$ and $H_2(z)$ are respectively $n_u$ and $\hat{n}_u$ and the degrees of the denominators of $H_1(z)$ and $H_2(z)$ are respectively $n_y$ and $\hat{n_y}$. Performing zero/pole cancellation does not change the difference between the degree of the numerator and the
degree of the denominator. Hence, we obtain that $n_y-n_u=\hat{n}_y-\hat{n}_u$ must hold. But since both $\AS_1$ and $\AS_2$ are minimal SARX realizations of the same input-output map, their dimensions must agree and hence $n_y+n_u=\hat{n}_y+\hat{n}_u$. It is easy to see that the only
  solution to the system of equations
 \[ \left\{\begin{array}{rcl}
             n_y-n_u&=&\hat{n}_y-\hat{n}_u, \\
             n_y+n_u&=&\hat{n}_y+\hat{n}_u,
           \end{array}\right.
 \]
 is $n_y=\hat{n}_y$ and $n_u=\hat{n}_u$.
 \end{proof} 

\ntext{As we have seen in the previous section, strong minimality implies
minimality. By \cite{Petreczky2012-Automatica,MP:HSCC2010}, strong minimality and hence minimality,  
can be checked algorithmically. Indeed, strong minimality of a SARX system $\AS$ means minimality of
the associated \LSS\ $\Sigma_{\AS}$. The latter can be checked by checking if the rank of each of the finite span-reachability matrix $\mathcal{R}(\Sigma_{\AS})$ of $\Sigma_{\AS}$ and the finite observability matrix $\mathcal{O}(\Sigma_{\AS})$ of $\Sigma_{\AS}$ considered in Theorem 2 ~\cite{Petreczky2012-Automatica} equals the dimension of $\Sigma_{\AS}$. We can also formulate sufficient conditions for minimality
which do not involve computing \SLSS.}
\begin{Theorem}[Sufficient conditions for (strong) minimality]
%\label{sarx2lss:lemma2}
\label{sarx2lss:lemma3}
%Consider the polynomial
 %\[ \chi_{q_1}(z)=z^{n_y}-\sum_{j=1}^{n_y} n_{q_1}^{j}z^{n_y-j}.\]
 Consider a SISO SARX system $\AS=\{h_q\}_{q \in Q}$ of type
 $(n_y,n_u)$.
  For all modes $q,\hat{q} \in Q$, define the polynomials
\begin{displaymath}
\chi_{q}(z)=z^{n_y}-\sum_{j=1}^{n_y} h_{q}^{j}z^{n_y-j}, \qquad
\upsilon_{q}(z)=\sum_{j=1}^{n_y} h_{q}^{j}z^{n_y-j},\qquad
\phi_{\hat{q},q}(z)=\sum_{j=1}^{n_u} h_{q}^{j+n_y}\psi_{\hat{q},q,n_u-j}(z),
\end{displaymath}  
with $\psi_{\hat{q},q,j}(z)$ defined recursively for $j=0,1,2\ldots,$ as follows: $\psi_{\hat{q},q,0}(z)=1$ and
 \begin{equation}
\label{sarx2lss:lemma2:eq1}
   \psi_{\hat{q},q,j+1}(z)=z\psi_{\hat{q},q,j}(z)+(h_{q}-h_{\hat{q}})d_j
  \end{equation}
where the vectors $d_j  \in \mathbb{R}^{n_y+n_u}$ are defined as follows:
$d_0=e_1$ and when $d_j=(d_{j,1},\ldots, d_{j,n_y},0,\ldots,0)^T$ with
$d_{j,1},\ldots, d_{j,n_y} \in \mathbb{R}$, $d_{j+1}=(h_qd_j, d_{j,1},\ldots, d_{j,n_y-1},0,\ldots,0)^T$. Then $\AS$ is strongly minimal, if the following conditions hold:
 \begin{description}
 \item{\textbf{(A)}}
  there exists discrete modes $q_0$ and $q_1$ such that
   the polynomials
 $\chi_{q_0}(z)$ and $\phi_{q_0,q_1}(z)$ are co-prime, and
 %then the \LSS\ $\Sigma_{\AS}$ is reachable.
 \item{\textbf{(B)}}
  there exists discrete modes $q_2$ and $q_3$, such that
 $\upsilon_{q_3}(z)$ and
 $\chi_{q_2}(z)$ are co-prime, $h_{q_2}^{n_u+n_y} \ne 0$ and 
  $h_{q_3}^{n_y} \ne \frac{h_{q_3}^{n_y+n_u}}{h_{q_2}^{n_y+n_u}}h_{q_2}^{n_y}$.
 %$\Sigma_{\AS}$ is observable.
%\item{\textbf{(C)}}
%  If the conditions \textbf{(A)} and \textbf{(B)} hold, then
%  $\AS$ is strongly minimal (and hence $\AS$ is minimal).
\end{description}
%% \begin{eqnarray*}
%%  \psi_{q_1}(z)&=& %z^{n_u}(
%%                  z^{n_y}-\sum_{j=1}^{n_y} n_{q_1}^{j} z^{n_y-j} \\
%%  \psi_{q_2}(z)&=& \sum_{i=1}^{n_y} n_{q_2}^{i} z^{n_y-i}
%%   %\sum_{i=1}^{n_u} n_{q_2}^{i+n_y} \upsilon_i(z)\phi(z) \\
%%  %\phi(z) &=& z^{n_y}-\sum_{i=1}^{n_y} n_{q_1}^{i} z^{n_y-i} \\
%%  %\upsilon_{1}(z)&=& \frac{1}{n_{q_1}^{n_u+n_y}} z^{n_u-1} \\
%%  %\upsilon_{i}(z)&=& \frac{1}{n_{q_1}^{n_u+n_y}}
%%   %(z^{n_u-i}-\sum_{j=1}^{i-1} \upsilon_{j}(z)n_{q_1}^{n_y+n_u-i+j})
%% \end{eqnarray*}
\end{Theorem}
\begin{proof}
See Appendix~\ref{appxTheo2}.
\end{proof}
\ntext{
\begin{remark}
Theorem \ref{sarx2lss:lemma3} is analogous to the well-known result that if a SISO transfer has no zero-pole cancellation (i.e. its numerator and denominator are coprime) and its denominator is of degree $n$, then all its minimal realizations are of order $n$. Due to the presence of switching, the formulation of Theorem \ref{sarx2lss:lemma3} is more involved. In addition, Theorem \ref{sarx2lss:lemma3} does not imply the classical results, since condition \textbf{(B)} of the theorem is always false if there is only one discrete state.
\end{remark}
 }

In order to demonstrate the utility of Theorem \ref{sarx2lss:lemma3}, we present the following examples.
 \begin{Example}\label{example:min3}
Let's apply Theorem \ref{sarx2lss:lemma3} to the SARX system $\AS$ \ntext{with two discrete modes $Q=\{1,2\}$ such that the ARX subsystem in mode $1$ is of the form $\y_{t}=8\y_{t-1}-15\y_{t-2}+\bu_{t-1}-3\bu_{t-2}$, and the ARX subsystem in mode $2$ is of the form $\y_{t}=\y_{t-1}+2\y_{t-2}+\bu_{t-1}+\bu_{t-2}$.} We obtain that $n_y=n_u=2$, $h_{1}=\big[\begin{matrix} 8 & -15 & 1 & -3 \end{matrix}\big]$ and $h_{2}=\big[\begin{matrix} 1 & 2 & 1 & 1 \end{matrix}\big]$. Hence, $h^{n_y+n_u}_{1}=-3 \ne 0$, and
    %\begin{eqnarray*}
    \(  \chi_{1}(z)=z^2-8z+15  \),
    \( \psi_{1,2,1}(z)= z-7 \)
     \( \upsilon_{2}(z)= z+2 \)
    \( \phi_{2,1}(z)= z-6 \).
  %\end{eqnarray*}
It is clear that the roots of $\chi_{1}(z)$ are $5$ and $3$ and hence $\upsilon_{2}(z)$ and $\chi_{1}(z)$ are co-prime and $\phi_{2,1}(z)$ and $\chi_{1}(z)$ are co-prime. Moreover, $h_{2}^{n_y}-\frac{h_{2}^{n_u+n_y}h_{1}^{n_y}}{h_{1}^{n_u+n_y}}=2-\frac{15}{3}=-3 \ne 0$. Hence, conditions \textbf{(A)} and \textbf{(B)} of Theorem \ref{sarx2lss:lemma3} hold and thus $\AS$ is (strongly) minimal.
 \end{Example}
\color{blue}
\begin{Example}
\label{Exple_Pbme:1}
 Consider the SARX system from Example \ref{Exple_Pbme}, and choose the parameter
 vector $\underline{\theta}=\begin{bmatrix} \theta_1 & \ldots & \theta_{20} \end{bmatrix}$ as
 \begin{equation}
 \label{Exple_Pbme:1:eq1} 
  \begin{split}
\underline{\theta}=& [  
 0.0046 , -0.0091, 0.0005,  -0.0019,  0.4881,  -0.9555, 0.0519,  
 -0.1973, -0.4881, 0.9555, -0.0519,  \\
 & 0.1973, 6.4616, -12.6262, 
 0.6924, -2.6043, -1.2564, 2.6133, -0.0989, 0.5625]
\end{split}
 \end{equation}
%We can apply Theorem \ref{sarx2lss:lemma3} to this example. 
It then follows that with this parameters, the SARX becomes
$\AS=\{n_q\}_{q \in Q}$, $Q=\{(i,j) \mid i,j=1,2\}$ with
$n_{(1,1)}=n_{(2,1)}=\begin{bmatrix} 8 & -15 & 1 & -3 \end{bmatrix}$ and 
$n_{(1,2)}=n_{(2,2)}=\begin{bmatrix} 1 & 2 & 1 & 1 \end{bmatrix}$.
Note that the parameter vectors are the same as in Example \ref{example:min3}
It follows that $\chi_{(1,1)}(z)=z^2-8z+15$, $\phi_{(1,1),(1,2)}(z)=z-6$.
 Hence, $\chi_{(1,1)}$ and $\phi_{(1,1),(1,2)}$ are co-prime and
condition \textbf{(A)} of Theorem \ref{sarx2lss:lemma3} holds for $q_0=(1,1)$ and $q_1=(1,2)$. 
Moreover, $\upsilon_{(1,2)}(z)=z+2$ and $\chi_{(1,1)}(z)=z^2-8z+15$ are also co-prime 
and $h_{(1,2)}^{n_y}=2$, $h_{(1,2)}^{n_y+n_u}=1$ and $ h_{(1,1)}^{n_y}=-15$, $h_{(1,1)}^{n_u+n_y}=-3$.
Hence,  $h_{(1,2)}^{n_y}-\frac{h_{(1,2)}^{n_u+n_y}h_{(1,1)}^{n_y}}{h_{(1,1)}^{n_u+n_y}}=2-\frac{15}{3}=-3 \ne 0$,
hence  condition \textbf{(B)} of Theorem \ref{sarx2lss:lemma3} holds for $q_2=(1,1)$ and $q_3=(1,2)$. 
That is, the SARX from Example \ref{Exple_Pbme} with the choice of parameters as in \eqref{Exple_Pbme:1:eq1} 
is strongly minimal.
\end{Example}
\color{black}

\subsection{Identifiability conditions for Switched ARX systems}
\label{ident:sarx}

In this section we study identifiability of SARX
systems. 
%\textcolor{blue}{and we show that strong minimality is sufficient for identifiability of SARX parametrizations. In order to state the result formally, we introduce the following definitions.}

%%%%%%%%%%%%%%%%%%%%%%%%%
\ntext{Theorem~\ref{sarx:main2.1}, which is one of the main results of the paper, describes the relationship between
 strong minimality and identifiability. More precisely, it show that strong minimality is sufficient for identifiability. In order to restrict attention to strongly minimal SARX systems, we introduce
the following terminology.}
\begin{Definition}[Minimality of SARX parametrizations]
  The parametrization $\MAPSARX$ 
  is called \emph{minimal} (resp. \emph{strongly minimal}), if
  for all $\theta \in \PARAMSET$, $\MAPSARX(\theta)$ is minimal
  (resp. strongly minimal).
 \end{Definition}
 \textcolor{blue}{If a SARX parametrization is strongly minimal, then the corresponding \SLSS\ parametrization will be minimal.  Hence, we can apply the conditions and algorithms provided by Petreczky {\it et al}~\cite{MP:HSCC2010} for analyzing the identifiability of the latter parametrization. By Corollary \ref{sarx2lss:col2} the identifiability of the latter parametrization is identifiability of the original SARX parametrization.}

In fact, for the SISO case (i.e. when $p=m=1$), we can 
derive even stronger results, by showing that 
\textit{minimality is sufficient for identifiability}.
 To this end, we need the following definition.
 %Another concept which we shall need is the one of
 %injective parametrizations.
\begin{Definition}{(Injective SARX parametrizations)}
 A SARX parametrization $\MAPSARX$ is said to be 
 \emph{injective}
 if $\MAPSARX$ is an injective map.
\end{Definition}
 An injective parametrization allows us to exclude 
 the situation
 where two different parameter values lead to the same
 SARX system. The ARX parametrization
 $\y_t=\theta^{2}\y_{t-1}+\bu_{t-1}$ with $\theta \in \mathbb{R}$ is not injective, since any 
 $\theta$ and $-\theta$ always lead to
 the same ARX system.

 \begin{Theorem}
 \label{sarx:main2.1}
  Assume that $p=m=1$.
  If a SISO SARX parametrization $\MAPSARX$ is injective and
 strongly minimal,
  then $\MAPSARX$ is identifiable.
 \end{Theorem}
\begin{proof}
See Appendix~\ref{appxTheo3}.
\end{proof} 
   % at least if some mild conditions hold.
%% \begin{pf}[Proof of Theorem \ref{sarx:main2.1}]
%%  From strong minimality of $\MAPSARX(\theta)$
%%  it follows that $\AS=\MAPSARX(\theta)$ is
%%  minimal and hence if $\AS=\{h_q\}_{q \in Q}$, then
%%  by Lemma \ref{element:lemma1} 
%%   for some $q \in Q$, $h_q^{n_u+n_y} \ne 0$. 
%%  Hence, by Theorem \ref{sarx:theo2},
%%  for any $\theta_1 \ne \theta_2 \in G$, the only 
%%  \LSS\ isomorphism between 
%%  $\Sigma_{\MAPSARX(\theta_1)}$ and 
%%  $\Sigma_{\MAPSARX(\theta_2)}$ is the identity. But the 
%%   latter implies that $\Sigma_{\MAPSARX(\theta_1)}=\Sigma_{\MAPSARX(\theta_2)}$ which implies that $\MAPSARX(\theta_1)=\MAPSARX(\theta_2)$. Since $\MAPSARX$ is an injective map, 
%%  $\theta_1=\theta_2$ which is a 
%%  contradiction. That is, there exists no isomorphism between
%%  $\Sigma_{\MAPSARX(\theta_1)}$ and $\Sigma_{\MAPSARX(\theta_2)}$.
%%  Hence, by \cite{MP:HSCC2010} the \LSS\ parametrization
%%  $\mathbf{\Pi}_{sw}:\PARAMSET \ni \theta \mapsto \Sigma_{\MAPSARX(\theta)}$ is identifiable.
%%\end{pf}
% Recall that strong minimality of a SARX system does not 
% imply minimality of its ARX subsystems. Similarly, 

\color{blue}
Theorem  \ref{sarx:main2.1} allows us to find an identifiable sub-parametrization of a SARX parametrization by checking finding a sub-parametrization which is strong minimal. 
One way to check strong minimality is by checking if the conditions of Theorem \ref{sarx2lss:lemma3} are satisfied.  This can easily be done, for
parameterizations in which the coefficients of the SARX systems depend on the parameters in a polynomial way. To this end, we introduce the following terminology. 
\color{black}
\begin{Definition}{(Polynomial parametrization)}  
  Let $K=(pn_y+mn_u)|Q|$. Then any SARX system of type $(n_y,n_u)$ can be identified with a point in $\mathbb{R}^{K}$, by identifying the system with its parameters $\{h_q\}_{q \in Q}$. Thus, $SARX(n_y,n_u,m,p,Q)$
 can be identified with the space $\mathbb{R}^{K}$.
 A parametrization $\MAPSARX$ is said to be \emph{polynomial}, if $\PARAMSET$ is an affine algebraic variety and $\MAPSARX$ is a polynomial map from $\PARAMSET$ to $SARX(n_y,n_u,m,p,Q)$. 
\end{Definition}
\color{blue}
Let $\MAPSARX$ be a polynomial parametrization. Below we present a procedure to find a subset $\hat{\Theta} \subseteq \Theta$ such that for each $\theta \in \hat{\Theta}$, the SARX
system $\MAPSARX(\theta)$ satisfies the conditions of Theorem \ref{sarx2lss:lemma3}, hence it is strongly minimal, and as a consequence the parametrization 
 $\PARAMSET|_{\hat{\Theta}}: \hat{\Theta} \ni \theta \mapsto \PARAMSET(\theta)$ is strongly minimal. To this end, we introduce the following notation.
We will use the standard notation and terminology from commutative algebra, see \cite{CoxBook}. In particular, we will need the notion of an ideal, generator of an ideal, Gr\"obner basi of an ideal, product of ideals from \cite{CoxBook}. 
 For each $\theta \in \PARAMSET$, if $\MAPSARX(\theta)=\{n_q(\theta)\}_{q \in Q}$, then denote by 
$\chi_q(\theta)(z)$, $\upsilon_{q}(\theta)(z)$ and $\phi_{q,\hat{q}}(\theta)(z)$ the polynomials $\chi_q(z),\upsilon_q(z), \phi_{q,\hat{q}}(z)$, $\hat{q},q \in Q$,
defined in Theorem \ref{sarx2lss:lemma3} for $n_q=n_{q}(\theta)$.  Then, since $\MAPSARX$ is polynomial, the dependence of $n_q(\theta)$ and the coefficients of $\chi_q(\theta)(z)$, $\upsilon_{q}(\theta)(z)$ and $\phi_{q,\hat{q}}(\theta)(z)$ on $\theta$ is polynomial. That is, 
there exist polynomials $n_q^{f,i} \in \mathbb{R}[X_1,\ldots,X_d]$, $i=1,\ldots,n_u+n_y$,  in variables $X_1,\ldots,X_d$, and polynomials
$\chi^f_q(X_1,\ldots,X_d,z),\upsilon^f_{q}(X_1,\ldots,X_d,z)$, $\phi_{q,\hat{q}}^f(X_1,\ldots,X_d,z)$
in variables $X_1,\ldots,X_d,z$ such that $n_q^i(\theta)=n_q^{f,i}(\theta)$, where $n_q^i(\theta)$ denotes the $i$th components of $n_q(\theta)$,  and
$i=1,\ldots,n_y+n_u$, $\chi_q(\theta)(z)=\chi^f_q(\theta,z)$, $\upsilon_{q}(\theta)=\upsilon^f_{q}(\theta,z)$ and 
$\phi_{q,\hat{q}}(\theta)(z)=\phi_{q,\hat{q}}^f(\theta,z)$, $\hat{q},q \in Q$. 
%Let $I \subseteq \mathbb{R}[X_1,\ldots,X_d]$ be the ideal of the variety $\PARAMSET$, i.e., $\PARAMSET=\{\theta \in \mathbb{R}^d \mid \forall P \in I: P(\theta)=0\}$. It is well-known \cite{} that $I$ is finitely generated by some polynomials $P_1,\ldots,P_K$, i.e., $\PARAMSET=\{\theta \in \mathbb{R}^d \mid P_i(x)=0, i=1,\ldots,K\}$. 
In order to apply Theorem \ref{sarx2lss:lemma3}, it is necessary to have a sufficient conditions for co-primeness of two polynomials in $z$, coefficients of which are polynomial functions of $\theta$. 
To this end,  assume that $Q_i(X_1,\ldots,X_d,z)$, $i=1,2$ are two polynomials. Consider the ideal $I(Q_1,Q_2)$ generated by the polynomials $Q_1,Q_2$ and consider
the ideal $J(Q_1,Q_2)=I(Q_1,Q_2) \cap \mathbb{R}[X_1,\ldots,X_d]$. The ideal $J(Q_1,Q_2)$ is finitely generated, and the set of its generators $S$ can be 
calculated from the polynomial $Q_1,Q_2$ using standard algorithms from compute algebra, see \cite{CoxBook} and the toolbox \cite{Sage}. 
\begin{Lemma}
\label{lemma:proc:ident0}
 If there exist $g \in S$ such that $g(\theta) \ne 0$, then the univariate polynomials $Q_i(\theta,z) \in \mathbb{R}[z]$, $i=1,2$ are co-prime. 
\end{Lemma}
\begin{proof}[Proof of Lemma \ref{lemma:proc:ident0}]
Indeed,  since
$g_i \in J(Q_1,Q_2) \subseteq I(Q_1,Q_2)$,
 $g(X_1,\ldots,X_d)=Q_1(X_1,\ldots,X_d,z)\alpha(X_1,\ldots,X_d,z)+Q_1(X_1,\ldots,X_d,z)\beta(X_1,\ldots,X_d,z)$ for some polynomials 
$\alpha,\beta$ in $R[X_1,\ldots,X_d]$. In particular, $g(\theta)=Q_1(\theta,z)\alpha(\theta,z)+Q_2(\theta,z)\beta(\theta,z)$ and since $g(\theta) \ne 0$,
$Q_1(\theta,z)\frac{\alpha(\theta,z)}{g(\theta)}+Q_2(\theta,z)\frac{\alpha(\theta,z)}{g(\theta)}=1$, which by Bezout's identity implies that $Q_1(\theta,z)$, 
$Q_2(\theta,z)$ are co-prime. 
\end{proof}
Lemma \ref{lemma:proc:ident0} implies that for any $\theta \in \{\theta \in \Theta \mid \exists P \in S: P(\theta) \ne 0\}$, the polynomials $Q_i(\theta,z) \in \mathbb{R}[z]$, $i=1,2$ are co-prime.
We can apply Lemma \ref{lemma:proc:ident0} to find $\hat{\Theta} \subseteq \Theta$ such that for any $\theta \in \hat{\Theta}$, the SARX system $\MAPSARX(\theta)$ satisfies 
conditions \textbf{(A)} and \textbf{(B)} of Theorem \ref{sarx2lss:lemma3}. More precisely, we propose the following 
algorithm for finding a strongly minimal sub-parametrization of a parametrization. 
\begin{procedure}[Identifiable polynomial parametrization]
\label{proc:ident}
\begin{enumerate}
\item
 For each $q,\hat{q} \in Q$, consider the ideal $I(\chi^f_q,\phi_{q,\hat{q}}^f)$ generated by the polynomials $\chi^f_q,\phi_{q,\hat{q}}^f$
  and calculate the Gr\"obner basis $S_{A,q,\hat{q}} \subseteq \mathbb{R}[X_1,\ldots,X_d]$ of the ideal $I(\chi^f_q,\phi_{q,\hat{q}}^f) \cap \mathbb{R}[X_1,\ldots,X_d]$
 using standard algorithms \cite{CoxBook}, implemented, for example, in the toolbox \cite{Sage}. 
 \item For each $q,\hat{q} \in Q$, consider $I(\chi_q^f, \upsilon_{\hat{q}}^f)$ generated by the polynomials $\chi^f_q,\upsilon_{\hat{q}}^f$
       and calculate the Gr\"obner basis  $S^{'}_{B,q,\hat{q}} \subseteq \mathbb{R}[X_1,\ldots,X_d]$ of the ideal 
       $I(\chi^f_q,\upsilon_{\hat{q}}^f) \cap \mathbb{R}[X_1,\ldots,X_d]$, 
       using standard algorithms \cite{CoxBook}, for an implementation see \cite{Sage}. 
       Let $S_{B,q,\hat{q}}=\{P\cdot Q \mid P \in S_{B,q,\hat{q}}^{'}\}$, where $Q=n_{q}^{f,n_u+n_y}(n_{\hat{q}}^{f,n_y}n_{q}^{f,n_u+n_y}-n_{q}^{f,n_y}n_{\hat{q}}^{n_y+n_u})$.

\item  Let $I_A$ be the ideal generated by $\bigcup_{q,\hat{q} \in Q} S_{A,q,\hat{q}}$,  and
       let $I_B$ be the ideal generated by $\bigcup_{q,\hat{q} \in Q} S_{B,q,\hat{q}}$ and let 
       $S$ be the Gr\"obner basis of the ideal $I_A I_B=\{P \mid P_1 \in I_A, P_2 \in I_B$. Note that $S$ can be computed from the Gr\"obner basis of $I_A$ and $I_B$,
       which, in turn, can easily be computed from the finite sets $\bigcup_{q,\hat{q} \in Q} S_{A,q,\hat{q}}$ and $\bigcup_{q,\hat{q} \in Q} S_{B,q,\hat{q}}$ using a standard
       algorithm for computing Gr\"obner basis from a generator set of an ideal \cite{CoxBook,Sage}. 
       Define the parametrization: $\MAPSARX|_{\hat{\Theta}}: \hat{\Theta} \ni \theta \mapsto \MAPSARX(\theta)$, where 
       \[ \hat{\Theta}=\{ \theta \in \Theta \mid \exists P \in S: P(\theta) \ne 0\}.  \]
\end{enumerate}
\end{procedure}
Procedure \ref{proc:ident} was implemented, the code is avaliable at \cite{ArxivePaper}. 
\begin{Lemma}
\label{lemma:proc:ident}
  The parametrization $\MAPSARX|_{\hat{\Theta}}$ calculated by Procedure \ref{lemma:proc:ident} is a  strongly minimal and hence it is identifiable
\end{Lemma}
\begin{proof}[Proof of Lemma \ref{lemma:proc:ident}]
 Assume that $\theta \in \hat{\Theta}$ and let $P \in S$ be such that $P(\theta) \ne 0$. Then, since $P \in I_AI_B$, $P=P_1P_2$ for some $P_1 \in I_A$ and $P_2 \in I_B$, and since $P(\theta) \ne 0$, $P_1(\theta) \ne 0$ and
 $P_2(\theta) \ne 0$. Since $P_1 \in I_A$ and $P_1(\theta) \ne 0$,  it then there mus exist a polynomial $\hat{P}_1$ in the generator set 
$\bigcup_{q,\hat{q} \in Q} S_{A,q,\hat{q}}$ of $I_A$ such that $\hat{P}_1(\theta) \ne 0$.  In particular, $\hat{P}_1 \in S_{A,q,\hat{q}}$ for some
 $q,\hat{q} \in Q$. 
 By applying Lemma \ref{lemma:proc:ident0} to $\chi^f_q,\phi_{q,\hat{q}}^f$ it follows that $\chi_q(\theta)(z)=\chi_q^f(\theta,z)$ and $\phi_{q,\hat{q}}^f(\theta,z)=\phi_{q,\hat{q}}(\theta)(z)$ are co-prime, hence
 for $q_1=q$, $q_2=\hat{q}$, condition \textbf{(A)} of Theorem \ref{sarx2lss:lemma3} holds.
 Similarly, since $P_2 \in I_B$ and $P_2(\theta) \ne 0$,
 it follows that  there exists a polynomial $\hat{P}_2$ such that $\hat{P}_2(\theta) \ne 0$ and for some $q_2,q_3 \in Q$, $\hat{P}_2 \in S_{B,q_2,q_3}$. The latter means that 
 $\hat{P}_2=P_3 \cdot n_{q_2}^{f,n_u+n_y}(n_{q_3}^{f,n_y}n_{q_2}^{f,n_u+n_y}-n_{q_2}^{f,n_y}n_{q_3}^{n_y+n_u})$ for some $P_3 \in S^{'}_{B,q_2,q_3}$. From $\hat{P}_2(\theta) \ne 0$ it then follows that
 $P_3(\theta) \ne 0$ and $n_{q_3}^{f,n_y}(\theta)n_{q_2}^{f,n_u+n_y}(\theta)-n_{q_2}^{f,n_y}(\theta)n_{q_3}^{n_y+n_u}(\theta) \ne 0$ and $n_{q_2}^{f,n_u+n_y}(\theta) \ne 0$. From Lemma \ref{lemma:proc:ident0}
 it follows that $P_3(\theta) \ne 0$ implies that $\upsilon_{q_3}^f(\theta)(z)$ and $\chi_{q_2}^f(\theta)(z)$ are co-prime, hence condition \textbf{(B)} of  Theorem \ref{sarx2lss:lemma3} holds.
 That is, $\PARAMSET(\theta)$ is strongly minimal. 
\end{proof}
\begin{Example}
\label{Exple_Pbme:2}
Consider the parametrization from Example \ref{Exple_Pbme} and let us apply Procedure \ref{proc:ident} to it.
We reparamaterize this parametrization as follows: define $\phi: \mathbb{R}^2 \ni (\zeta_1,\zeta_2)^T \mapsto (\zeta_1,\zeta_1+\zeta_2,0,0,\ldots,0,1)^T \in \mathbb{R}^{20}$
and define the parametrizaton $\bar{\MAPSARX}: \mathbb{R}^2 \ni (\zeta_1,\zeta_2)^T \mapsto \MAPSARX(\phi(\zeta_1,\zeta_2))$, where $\MAPSARX$ is the parametrization
from Example \ref{Exple_Pbme}. It then follows that 
for $\underline{\zeta}=(\zeta_1,\zeta_2)$  are of the form
$\hat{\MAPSARX}(\underline{\zeta})=\{n_q(\underline{\zeta})\}_{q \in Q}$, $Q=\{(1,1),(1,2),(2,1),(2,2)\}$,
\[
  \begin{split}
    n_{(1,2)}(\underline{\zeta})= \begin{bmatrix} 0.001\zeta_1 & 0.001\zeta_1 + 0.001\zeta_2 &  0.0 & 16.0 \end{bmatrix}^T, ~ 
    n_{(1, 1)}(\underline{\zeta})= \begin{bmatrix} 0.001\zeta_1 & 0.001\zeta_1 + 0.001\zeta_2 & 0.0 & 4.0 \end{bmatrix}^T \\
    n_{(2, 1)}(\underline{\zeta})= \begin{bmatrix} 0.001\zeta_1 & 0.001\zeta_1 + 0.001\zeta_2 & 0.0 & 4.0 \end{bmatrix}^T, ~
    n_{(2, 2)}(\underline{\zeta})= \begin{bmatrix} 0.001\zeta_1 & 0.001\zeta_1 + 0.001\zeta_2 & 0.0 & 16.0 \end{bmatrix}^T
  \end{split}
\]
and the polynomials $\chi_q(\underline{\zeta}), \upsilon_q(\underline{\zeta}), \phi_{q_1,q_2}(\underline{\zeta})$,
\[
   \begin{split}
    & \chi_{q}(\underline{\zeta})(z)=-0.001\zeta_1 z + z^2 - 0.001\zeta_1 - 0.001\zeta_2, ~ 
    \upsilon_{q}(\underline{\zeta})(z)=0.001\zeta_1 z  + 0.001\zeta_1 + 0.001\zeta_2, ~ q \in Q \\
    & \phi_{q_1,q_2}(\underline{\zeta})(z)=16z, ~ (q_1,q_2) \in Z_1, ~ 
    \phi_{q_1,q_2}(\underline{\zeta})(z)=4z, ~ (q_1,q_2) \in Z_2, \\
  &Z_1=\{((1,2),(2,2)),((1,1),(2,1)),((1,2),(2,2)),((2,1),(1,1)),((2,2),(1,2))\} \\
  & Z_2=\{((1,1),(1,2)),((1,1),(2,2)), ((1,2),(1,1)),((2,1),(1,2)), ((2,1),(2,2)),((2,2),(1,1))\}.
    %\phi_{(1,1),(1,2)}(\underline{\zeta})(z)=\phi_{(1,1),(2,2)}(\underline{\zeta})(z)=\phi_{(1,2),(2,2)}(\underline{\zeta})(z)=
   %\phi_{(2,1),(2,2)}(\underline{\zeta})(z)=\phi_{(2,1),(1,2)}(\underline{\zeta})(z)=\phi_{(2,1),(2,2)}(\underline{\zeta})(z)=
%\phi_{(2,2),(1,2)}(\underline{\zeta})(z)=16.0z \\
%    \phi_{(1,1),(2,1)}(\underline{\zeta})(z)=\phi_{(1,2),(1,1)}(\underline{\zeta})(z)=\phi_{(1,2),(2,1)}(\underline{\zeta})(z)=
%  \phi_{(2,1),(1,1)}(\underline{\zeta})(z)=\phi_{(2,2),(1,1)}(\underline{\zeta})(z)=\phi_{(2,2),(1,2)})(\underline{\zeta})(z)=4.0z
   \end{split}
\]
Then $S_{A,q,\hat{q}}$, $S_{B,q,\hat{q}}^{'}$ and $S_{B,q,\hat{q}}$, $q,\hat{q} \in Q$ can all be calculated using \cite{Sage}. 
In this case $S_{A,q,\hat{q}}=\{ \zeta_1+\zeta_2 \}$, $q,\hat{q} \in Q$, $q \ne \hat{q}$ and 
$S_{B,q,\hat{q}}=\emptyset$, $(q,\hat{q}) \in Z_1$, 
$S_{B,q,\hat{q}}=\zeta_1^3 + 3\zeta_1^2\zeta_2 + 3\zeta_1\zeta_2^2 + \zeta_2^3$,
 $(q,\hat{q}) \in Z_2$.
It then follows that $I_A$ is generated by the Gr\"obner basis  $\{\zeta_1+\zeta_2\}$ and $I_B$ is generated by 
the Gr\"obner basis $\{\zeta_1^3 + 3\zeta_1^2\zeta_2 + 3\zeta_1\zeta_2^2 + \zeta_2^3\}$. Hence, 
$S=\{\zeta_1^4 + 4\zeta_1^3\zeta_2 + 6\zeta_1^2\zeta_2^2 + 4\zeta_1\zeta_2^3 + \zeta_2^4\}$
It then follows that 
\[ \hat{\Theta}=\{ \underline{z}=(\zeta_1,\zeta_2) \mid \zeta_1^4 + 4\zeta_1^3\zeta_2 + 6\zeta_1^2\zeta_2^2 + 4\zeta_1\zeta_2^3 + \zeta_2^4 \ne 0\}=
    \{ \underline{z}=(\zeta_1,\zeta_2) \mid \zeta_1 \ne -\zeta_2 \}
\]
and the parametrization
$\bar{\MAPSARX}|_{\hat{\Theta}}: \hat{\Theta} \ni \underline{\zeta} \mapsto \bar{\MAPSARX}(\underline{\zeta})$ is
strongly minimal and identifiable.

We can also apply Procedure \ref{proc:ident} to more complicated parametrizations, but the expressions for corresponding
polynomials and the Gr\"obner bases are more invoved. For example, 
define $\tilde{\phi}: \mathbb{R}^2 \ni (\zeta_1,\zeta_2)^T \mapsto (\zeta_1,\zeta_1+\zeta_2,\kappa_3,\kappa_4,\ldots,\kappa_{18},\kappa_{19}+1)^T \in \mathbb{R}^{20}$, $\kappa_{i}=1$, if $i$ is even, and $\kappa_i=(\zeta_1-\zeta_2)$ if $i$ is odd, 
and define the parametrizaton $\tilde{\MAPSARX}: \mathbb{R}^2 \ni (\zeta_1,\zeta_2)^T \mapsto \MAPSARX(\tilde{\phi}(\zeta_1,\zeta_2))$, where $\MAPSARX$ is the parametrization from Example \ref{Exple_Pbme}.  The expressions for the Gr\"obner basis 
$S_{A,q,\hat{q}}$, $S_{B,q,\hat{q}}$, $q, \hat{q} \in Q$ is lengthy. However, using the implementation of \cite{ArxivePaper}
of Procedure \ref{proc:ident}, we obtain that $I_A$ is generated by the polynomial $1$, i.e., $I_A=\mathbb{R}[X_1,X_2]$,
and $I_B$ is generated by $\{\zeta_1^2,\zeta_1\zeta_2,\zeta_2^2\}$ and hence $I_AI_B=I_B$ and thus
$S=\{\zeta_1^2,\zeta_1\zeta_2,\zeta_2^2\}$ and  
\[ \tilde{\Theta}=\{ \underline{z}=(\zeta_1,\zeta_2) \mid \zeta_1^2 \ne 0, \mbox{ or } \zeta_1\zeta_2  \ne 0 \mbox{ or } \zeta_2^2 \ne 0\}=
    \{ \underline{z}=(\zeta_1,\zeta_2) \mid \zeta_1 \ne 0 \mbox{ or }  \zeta_2 \ne 0 \}
\]
and the parametrization
$\tilde{\MAPSARX}|_{\tilde{\Theta}}: \tilde{\Theta} \ni \underline{\zeta} \mapsto \tilde{\MAPSARX}(\underline{\zeta})$ is
strongly minimal and identifiable. 
\end{Example}
\begin{remark}[Computional complexity]
  Procedure \ref{proc:ident} relies on computing Gr\"obner bases, and it is known that the computational complexity of the
  latter can be high. Hence, computational complexity of Procedure \ref{proc:ident} might be an issue for applications.
  However, even for linear systems, identifiability analysis relies on symbolic algorithms, in particular, on
  algorithms based on calculation of Gr\"obner bases, and there the same problem arises \cite{JvandenHof1998}. For this reason, 
  a detailed study of computational complexity of Procedure \ref{proc:ident} cannot be handled within this paper. 
\end{remark}
\color{black}

%\subsection{Identifiability of SARX systems}
\subsection{On the genericity of minimality and identifiability}
\label{main:generic}
  In the previous sections we have established that
  strong minimality is sufficient for minimality and that
  it is also sufficient for identifiability. However, we
  have also demonstrated that for some 
  minimal SARX systems, strong minimality
  does not hold. Hence, one may wonder how typical strong
  minimality is.

 Below we will show that strong minimality is a generic 
 property, i.e. it holds for almost all SARX systems, if
 $|Q|>1$. This also means that identifiability is a 
 generic property.  In other words, strong minimality
 occurs very frequently. In order to formalize these results, we need the following
terminology.
\begin{Definition}[Generic set]
 A subset $G$ of $\PARAMSET \subset \mathbb{R}^{d}$ is \emph{generic}, if $G$ is non-empty and
 there exists a non-zero polynomial $P(X_1,\ldots,X_d)$ in
 $d$ variables such that 
 $G=\{ \theta \in \PARAMSET \mid P(\theta)\ne 0\}$. 
 %$\theta \notiin G$, $ \theta \in \PARAMSET$, $P(\theta)=0$.
\end{Definition}
 That is, a generic subset of $\PARAMSET$ is a non-empty subset 
 whose complement in $\PARAMSET$ satisfies a polynomial
 equation.
\begin{Definition}[Generic identifiability and minimality of SARX parametrization]
 The parametrization $\MAPSARX$  is said to be
 \emph{generically identifiable} 
 if there exists a generic subset
 $G$ of $\PARAMSET$, such that the
 parametrization $\MAPSARX|_{G}:G \ni \theta \mapsto \MAPSARX(\theta)$ is identifiable. Similarly, $\MAPSARX$ is \emph{generically minimal} (respectively \emph{generically strongly minimal}), if there exists a generic subset $G$ of $\PARAMSET$,
 such that the parametrization $\MAPSARX_{|G}: G \ni \theta \mapsto \MAPSARX(\theta)$ is minimal (respectively strongly minimal).
\end{Definition}

 Intuitively, if a property is generic for a parametrization, then every member of the parametrization can be approximated with arbitrary accuracy by another member which has  this property. Another interpretation is that if we randomly generate parameters, then the property will hold for the obtained random parametrization with probability one.

\begin{Example}
 Consider the parametrization $\MAPSARX$ from 
  Example \ref{ident:example}. The set $G$ from
  Example \ref{ident:example} is generic. Hence,
  since the parametrization $\MAPSARX|_{G}$ is 
  strongly minimal and identifiable, the parametrization
  $\MAPSARX$ is generically strongly minimal,
  generically minimal, and generically identifiable.
 \end{Example}

 \begin{Theorem}[Generic minimality]
 \label{sarx:theo1}
  If $|Q|>1$, $\MAPSARX$ is a polynomial parametrization and
  $\MAPSARX$ contains a strongly minimal SARX system, \
 (i.e. for some $\theta \in \PARAMSET$, $\MAPSARX(\theta)$ is strongly minimal), 
  then $\MAPSARX$ is \emph{generically strongly minimal}.
 \end{Theorem}
\begin{proof}
See Appendix~\ref{appxTheo5}.
\end{proof}
Notice that
 Theorem \ref{sarx:main2.1} implies the following corollary result.
 \begin{Corollary}
 \label{sarx:col_main2}
  Consider the SISO case, i.e. $p=m=1$.
  If a SARX parametrization is injective, polynomial, and  generically
  strongly minimal, then it is generically identifiable.
 \end{Corollary}
 \begin{proof}%[Proof of Corollary \ref{sarx:col_main2}]
  If $\MAPSARX$ is generically strongly minimal, then there
  exists a generic set $G \subseteq \PARAMSET$ such
  that the parametrization 
  $\MAPSARX|_{G}:G \ni \theta \mapsto \MAPSARX(\theta)$
  is strongly minimal. Hence, by Theorem \ref{sarx:main2.1},
  $\MAPSARX|_{G}$ is identifiable. This means that 
  $\MAPSARX$ is generically strongly identifiable.
 \end{proof}
%%%
  Corollary \ref{sarx:col_main2} and Theorem \ref{sarx:theo1} yield the following result.
 \begin{Corollary}
 \label{sarx:main2}
  Assume that $p=m=1$.
  If a SISO SARX parametrization is polynomial and it
  contains a strongly minimal element, 
  then it is generically identifiable.
 \end{Corollary}
%%%
 %\begin{proof}%[Proof of Corollary \ref{sarx:main2}]
  %If $\MAPSARX$ contains a strongly minimal element,
  %then by Theorem \ref{sarx:theo1}, $\MAPSARX$ is
  %generically strongly minimal. The rest of the proof is a direct consequence of 
  %Corollary \ref{sarx:col_main2}.
 %\end{proof}
%%%
The trivial SISO SARX parametrization $\mathbf{\Pi}_{\text{\scriptsize triv}}$
  is the SARX parametrization defined as follows: $\PARAMSET=\mathbb{R}^{|Q|(n_u + n_y)}$ and $\mathbf{\Pi}_{\text{\scriptsize triv}}$ is the identity map.
From Corollaries \ref{sarx:col_main2} and \ref{sarx:main2}, we obtain that
%%%%
 \begin{Corollary}
 \label{sarx:main:col2}
  The trivial parametrization is generically minimal and in the SISO case, it is generically identifiable.
 \end{Corollary}
%%%
 \begin{proof}%[Proof of Corollary \ref{sarx:main:col2}]
  By Remark \ref{example:min2}, there exists a strongly minimal
  SARX system, i.e. $\mathbf{\Pi}_{\text{\scriptsize triv}}$ contains a strongly minimal element. Moreover, $\mathbf{\Pi}_{\text{\scriptsize triv}}$ is clearly injective and polynomial. We can therefore apply Corollary \ref{sarx:main2} to conclude.
\end{proof}
\color{blue}
\begin{Example}
\label{Exple_Pbme:3}
 From Example \ref{Exple_Pbme:1} it follows that the parametrization defined in Example \ref{Exple_Pbme} contains a strongly minimal element,
 hence it is generically strongly minimal and generically identifiable.
\end{Example}
\color{black}

\section{Concluding remarks}
\label{sect:concl}
In this paper we have studied  minimality and identifiability of linear SARX systems. Formal definitions of these two concepts have been introduced and discussed with respect to their standard characterizations for ARX systems. Sufficient and necessary conditions have been derived for minimality and identifiability of SARX systems. In particular, it has been shown that minimal SARX parametrizations are also identifiable. 

%\section*{Acknowledgments}
%This work is supported in part by a FNADT subvention in the
%framework of the SUNRISE project.

\subsection*{ORCID}

\appendix

\expandafter\def\expandafter\ref\expandafter#\expandafter1\expandafter{\expandafter\lowercase\expandafter{\ref{#1}}}

\section{On minimality and indentifiability of \SLSS}\label{AppdxLSS}

%This part recalled the formal definition of minimality and indentifiability of \SLSS\ as well as some important characterization results of these properties. Throughout this section, $\Sigma=\SwitchSysLin$ denotes a \SLSS\  of the form (\ref{lin_switch0}). Note that as for the notation introduced for SARX systems, in this notation we did not include explicitly $m$ and $p$, as these integers are fixed in the paper. Furthermore, a sequence $w=(q_0,u_0)\cdots (q_t,u_t) \in \HYBINP^{+}$ describes the scenario, when discrete mode $q_i$ and continuous input $u_i$ are fed to $\Sigma$ at time $i$, for $i=0,\ldots,t$.}
%\subsection{\SLSS\, input-output maps realization and minimality}
    %\begin{Definition}[Input-output maps of \SLSS]
   
 %  \end{Definition}
%%%%%%%%%%%%%%%%%%%%%%%%%%%%%%%%%%%%%%%%%%%%%%%%%%%%%%%%%%%%%%%%%%%%%%%%%%%%%%%%
%   \begin{Definition}[Dimension of \SLSS]
%   \label{switch_sys:dim:def}
   %\end{Definition}
 %%\vspace{-10pt}
 %Note that the number of discrete states is fixed, and
 %hence not included into the definition of dimension.
%\begin{Definition}{(Minimality of \BSLSS)}
%\label{def:minimality_LSS}
%\end{Definition}
 %%\vspace{-10pt}
\color{blue}
%\subsection{Structural identifiability and minimality of \SLSS}
\label{sect:ident}
\label{sect:Identifiability}
 In this section we recall from \cite{MP:HSCC2010} the notion of identifiability for
 \SLSS, and its relationship with minimality. In addition,
  we recall from \cite{MP:HSCC2010} sufficient and necessary conditions for
  identifiability of \SLSS.
%%%%%
% \paragraph{Structural identifiability}
 %We define the notion of structural
 %identifiability of parametrizations of \SLSS. 
 We start with defining 
 the notion of parametrization of \SLSS. To this end, we need the
 following notation.
%%%%%%%%%%%%%%%%%%
 \begin{Notation}
  Denote by $\Sigma(n,m,p,Q)$ the set of all \SLSS\ with
  state-space dimension $n$, input space $\mathbb{R}^{m}$,
  output space $\mathbb{R}^{p}$, and set of discrete modes $Q$.
 \end{Notation}
%%%%
 \begin{Definition}[Parametrization of \SLSS]
 Assume that $\PARAMSET \subseteq \mathbb{R}^{d}$ is the set of parameters.
 A \emph{parametrization of \BSLSS belonging to $\Sigma(n,m,p,Q)$} is a map
 $ \MAPLSS:\PARAMSET \rightarrow \Sigma(n,m,p,Q)$. 
  For each $\theta \in \PARAMSET$, 
  we denote $\MAPLSS(\theta)$ by 
  \[ \Sigma(\theta)=\PARAMLinSwitch{\theta}{{}}. \]
 \end{Definition}
%%%%
 Next, we define structural identifiability of parametrizations.
 %This notion is central to the study of identifiability of
 %\SLSS.
%%%%
 \begin{Definition}[Structural identifiability of \SLSS\, parametrizations]
  A parametrization
  $\MAPLSS:\PARAMSET \rightarrow \Sigma(n,m,p,Q)$ is 
  \emph{structurally identifiable}, if 
  for any two distinct parameters $\theta_1 \ne \theta_2$,
  the input-output maps of the corresponding
  \BSLSS $\MAPLSS(\theta_1)=\Sigma(\theta_1)$ and 
  $\MAPLSS(\theta_2)=\Sigma(\theta_2)$
  are different, i.e. 
  \( y_{\Sigma(\theta_1)} \ne y_{\Sigma(\theta_2)} \).
   %\[ \PARAMSET \ni \theta \mapsto y_{\Sigma(\theta)}: \HYBINP^{+} \rightarrow \mathbb{R}^{p} \mbox{ where } \Sigma(\theta)=\MAPLSS(\theta)  \]
   %is injective.
   %Recall that $y_{\Sigma(\theta)}$ denotes the input-output map
   %of the \BLSS $\Sigma(\theta)$.
 \end{Definition}
%%%%\vspace{-10pt}
  The condition $y_{\Sigma(\theta_1)}\ne y_{\Sigma(\theta_2)}$ means
  that there exists a sequence of inputs and discrete modes
  $w \in \HYBINP^{+}$, such that
  $y_{\Sigma(\theta_1)}(w) \ne y_{\Sigma(\theta_2)}(w)$.
  In other words, a parametrization is structurally identifiable, if for every two distinct parameters there exists an input and
  a switching signal, such that the corresponding outputs are different.
  This means that every parameter can be uniquely reconstructed from the input-output map of the corresponding LSS.
%%%%

% \paragraph{Identifiability and minimality}

It is an intuitive fact that minimality is somehow a necessary condition for structural identifiability \cite{MP:HSCC2010}. If we allow non-minimal parametrizations, then either the parametrization is not identifiable, or all the parameters occur in the minimal part of the systems, and hence we can replace the parametrization by a minimal one. For this reason, \emph{we will restrict attention to minimal \BSLSS when studying identifiability}. In turn, structural minimality of parametrizations allow a simple characterization of identifiability, due to the fact that minimal \BSLSS are unique up to isomorphism.

 \begin{Definition}[Structural minimality of \SLSS\, parametrization]
  The parametrization $\MAPLSS$ is called structurally minimal, if
 for any parameter value $\theta \in \PARAMSET$, $\Sigma(\theta)$ is a minimal \BLSS realization of its input-output map $y_{\Sigma(\theta)}$. 
 \end{Definition}
%%%%\vspace{-10pt}
  Hence, by Petreczky {\it et al}~\cite[Theorem 1]{Petreczky2012-Automatica}, 
  $\MAPLSS$ is structurally minimal if and only if for every parameter $\theta \in \Theta$, $\Sigma(\theta)$ is \WR\ and observable. Since the latter concepts admit rank characterizations, structural minimality is a property that can be checked algorithmically.

%\paragraph{Characterization of structural identifiability}

Theorem~\ref{theo:structural_ident} below recalls a necessary and sufficient condition for structural identifiability of a \emph{structurally minimal} parametrization established by Petreczky {\it et al}~\cite{MP:HSCC2010}.

  \begin{Theorem}[Identifiability of structural minimal parametrizations]
  \label{theo:structural_ident}
   A structurally minimal parametrization $\MAPLSS$ is structurally identifiable, if and only if for any two distinct parameter values $\theta_1, \theta_2 \in \PARAMSET$, $\theta_1 \ne \theta_2$, there exists no \BLSS isomorphism $S: \Sigma(\theta_1) \rightarrow \Sigma(\theta_2)$.
  \end{Theorem}
  
 The following important corollary which is an immediate consequence of the Theorem~\ref{theo:structural_ident} can be useful for checking identifiability of parametrizations.
  %%\vspace{-20pt}
 \begin{Corollary}
  \label{theo:structural_ident:col}
  Assume that $\MAPLSS$ is a structurally minimal parametrization, 
  and for each two parameter values
  $\theta_1,\theta_2 \in \PARAMSET$,
  $\Sigma(\theta_1)=\Sigma(\theta_2)$
  implies that
  $\theta_1=\theta_2$. 
  Here, equality of two systems means equality of the matrices of 
  the linear subsystems for each discrete state $q \in Q$ and 
  equality of the initial state.
  Then $\MAPLSS$ is structurally identifiable
  if and only if the assumption that 
  $S:\Sigma(\theta_1) \rightarrow \Sigma(\theta_2)$ is an \BLSS isomorphism
  implies that $S$ is the identity matrix.
 \end{Corollary}
\color{black}

\section{Proof of Theorem~\ref{sarx2lss:lemma3}}\label{appxTheo2}
  For the proof of Theorem \ref{sarx2lss:lemma3}, we will 
  need a number of auxiliary results.
  Below, we consider $\AS=\{h_q\}_{q \in Q}$.
   We denote by $A_q$ the corresponding matrix of the 
  \LSS\ $\Sigma_{\AS}$.
   We will denote by $e_i$ the $i$th standard basis vector of $\mathbb{R}^{n_y+n_u}$.
 \begin{Lemma}
 \label{sarx2lss:lemma2:lemma2.1}
  Let $X_1=\mathrm{Span}\{e_1,\ldots,e_{n_y}\}$. It then follows
  that for any $q \in Q$,
  \begin{enumerate}
  \item
  \label{sarx2lss:lemma2:lemma2.1:part1}
    $A_{q}e_j=h_{q}^{j}e_1+e_{j+1}$ for all
    $j=1,\ldots, n_y+n_u$, $j \ne n_y$ and
    $A_{q}e_{n_y}=h_{q}^{n_y}e_1$.
  \item  
 \label{sarx2lss:lemma2:lemma2.1:part2}
  The space $X_1$ is $A_{q}$ invariant and
  $\mathrm{Span}\{A_{q}^{i}e_1 \mid i=0,\ldots,n_y-1\}=X_1$
  \item
 \label{sarx2lss:lemma2:lemma2.1:part3}
  $A_{q}^{j}e_{n_y+1} \in X_1+e_{n_y+j+1}$, $j=0,\ldots,n_u-1$,
  and $A_{q}^{n_u}e_{n_y+1} \in X_1$.
 \item
 \label{sarx2lss:lemma2:lemma2.1:part4}
  For any $\hat{q} \in Q$, $\psi_{\hat{q},q,j}(A_{\hat{q}})e_1=A_{q}^je_1$,
  where $\psi_{\hat{q},q,0}(z)=1$ and
 \[ \psi_{\hat{q},q,j+1}(z)=z\psi_{\hat{q},q,j}(z)
    +(h_{q}-h_{\hat{q}})\psi_{\hat{q},q,j}(A_{\hat{q}})e_1, 
  \]
 where the vectors $d_j \in \mathbb{R}^{n_y+n_u}$ 
 are defined as follows:
$d_0=e_1$ and if $d_j=(d_{j,1},\ldots, d_{j,n_y},0,\ldots,0)^T$,
$d_{j,1},\ldots, d_{j,n_y} \in \mathbb{R}$, 
then 
\[ d_{j+1}=(h_qd_j, d_{j,1},\ldots, d_{j,n_y-1},0,\ldots,0)^T \]
 \end{enumerate}
 \end{Lemma}
 \begin{proof}%[Proof of Lemma \ref{sarx2lss:lemma2:lemma2.1}]
   Part \ref{sarx2lss:lemma2:lemma2.1:part1} follows by a 
   simple computation. 
   Part \ref{sarx2lss:lemma2:lemma2.1:part2} follows
   from Part \ref{sarx2lss:lemma2:lemma2.1:part1} by taking
   into account that $A_{q}e_{n_y} \in \mathrm{Span}\{e_1\}$.
   Part \ref{sarx2lss:lemma2:lemma2.1:part3} follows from
   the definition of $A_q$ by induction. Indeed, $A_{q}e_{n_y+1}=h_{q}^{n_y+1}e_1+e_{n_y+2} \in X_1 + e_{n_y+2}$ and if $A_{q}^{j}e_{n_y+1} \in X_1+e_{j+1}$, then 
  $A_{q}^{j+1}e_{n_y+1} \in X_1 + A_{q}e_{j+1} \subseteq X_1+h_q^{n_y+j+2}e_1+e_{j+2} \subseteq X_1+e_{j+2}$.
  %Notice that this implies that $A_{q}^{j}e_{n_y+1}$,
  %$j=0,\ldots,n_u-1$ are linearly independent.
  %Notice that $A_{q}e_{n_y+n_u}=h_{q}^{n_y+n_u}e_1$. Since
  %$h_{q}^{n_y+n_u} \ne 0$, we can conclude that
 % $e_1  \in \mathrm{Span}\{A_{q}^{j}e_{n_y+1} \mid j=0,\ldots,n_u\}$.

  Finally, $\psi_{\hat{q},q,j}(A_{\hat{q}})e_1=A_{q}^{j}e_1$ we will
  prove by induction. For $j=0$, the equality is trivial.
  Notice that 
  \( A_{q}e_i=h_{q}^{i}e_1+e_{i+1}=A_{\hat{q}}e_i+(h_{q}^{i}-h^{i}_{\hat{q}})e_1 \) for all $i=1,\ldots,n_y-1$, and
   \( A_{q}e_{n_y}=h_{q}^{n_y}e_1=A_{\hat{q}}e_{n_y}+(h_{q}^{n_y}-h_{\hat{q}}^{n_y})e_1 \). Hence, for any $x=\sum_{i=1}^{n_y} x_ie_i$, 
  \( A_{q}x=A_{\hat{q}}x+\sum_{i=1}^{n_y} x_i(h_{q}^{i}-h_{\hat{q}}^i)e_1=
     A_{q_1}x+((h_{q}-h_{\hat{q}})x)e_1
  \). Hence, if $\psi_{\hat{q},q,j}(A_{\hat{q}})e_1=A_{q}^je_1$ holds, then
  \begin{equation}
  \label{sarx2lss:lemma2:lemma2.1:part4:eq1}
     A_{q}^{j+1}e_1=A_{q}\psi_{\hat{q},q,j}(A_{\hat{q}})e_1=
     A_{\hat{q}}\psi_{\hat{q},q,j}(A_{\hat{q}})e_1+
      ((h_q-h_{\hat{q}})\psi_{\hat{q},q,j}(A_{\hat{q}})e_1)e_1.
  \end{equation}
  Finally,  notice that $d_j=A_{q}^{j}e_j$ for all $j=0,\ldots,n_y$.
  Indeed, $d_0=e_1$ and if $d_j=\sum_{i=1}^{n_y} d_{j,i}e_i$,
  then $A_{q}d_j=(\sum_{i=1}^{n_y} d_{j,i}h_{q}^i)e_1+\sum_{i=2}^{n_y} d_{j,i-1}e_i=d_{j+1}$.
  Hence, by replacing $\psi_{\hat{q},q,j}(A_{\hat{q}})e_1=A_{q}^{j}e_1$ by $d_j$ in
  \eqref{sarx2lss:lemma2:lemma2.1:part4:eq1}, we obtain that
  \[ A_{q}^{j+1}e_1=\psi_{\hat{q},q,j+1}(A_{\hat{q}})e_1. \]
  Hence, by induction we get the last statement of the lemma.
%%  By induction it is easy to see that 
%%  $A_{q}^{j}e_1 \in \mathrm{Span}\{e_1,\ldots,e_{j-1}\}+e_{j}$,
%%  $j=0,\ldots,n_{y}-1$. Hence, 
%%  $A_{q}^{j}e_1$, $j=0,\ldots,n_{y}-1$ are linearly independent
%%  and together they span $X_1=\mathrm{Span}\{e_1,\ldots,e_{n_y}\}$.
%%  Combining all this we get that the vector
%%  $A_{q}^{n_u+j}e_{n_y+1}$, $j=0,\ldots,n_y-1$ span
%%  $X_1$. Since $A_{q}^{j}e_{n_y+1} \in X_1 + e_{n_y+j+1}$ for
%%  $j=0,\ldots,n_{u}-1$, we then obtain that
%%  the vectors $A_{q}^{j}e_{n_y+1}$, $j=0,\ldots,n_y+n_u-1$ span
%%  the whole $\mathbb{R}^{n_y+n_u}$. Since 
%%  $B_q=e_{n_y+1}$, this means that $(A_q,B_q)$ is a controllable
%%  pair, from which by Theorem \ref{sect:real:lemma1}
%%  it follows that
%%  $\Sigma_{\AS}$ is reachable.
 \end{proof}
  %Below we assume that some $q \in Q$ is such that $h_{q}^{n_y+n_u} \ne 0$.
  \begin{Lemma}
  \label{sarx2lss:lemma3:lemma1}
   If $h_{q}^{n_y+n_u} \ne 0$, then
   $\{e_{n_y}^TA_{q}^{j} \mid j=0,\ldots,n_y+n_u-1 \}$ spans
   $\mathbb{R}^{1 \times (n_y+n_u)}$. Moreover, $e_i^T=e_{n_y}^TA^{n_y-i}_{q}$, $i=1,\ldots,n_y$,
   and $e_{n_y+j}^T=e_{n_y}^T\chi_q(A_q)\gamma_{j,q}(A_q)$, $j=1,\ldots,n_u$,
    where $\chi_{q}(z)=z^{n_y}-\sum_{j=1}^{n_y} h_q^{j}z^{n_y-j}$ and  the polynomial $\upsilon_j(z)$, $j=1,\ldots,n_u$  
   is defined recursively as follows:
\begin{displaymath}
\gamma_{1,q}(z)=\frac{1}{h_{q}^{n_u+n_y}} z^{n_u-1} 
\end{displaymath}
and
\begin{displaymath}
\gamma_{i,q}(z)= \frac{1}{h_{q}^{n_u+n_y}}
    (z^{n_u-i}-\sum_{j=1}^{i-1} \gamma_{j,q}(z)h_{q}^{n_y+n_u-i+j}).
\end{displaymath}
%   \begin{eqnarray*}
%   \gamma_{1,q}(z)&=& \frac{1}{h_{q}^{n_u+n_y}} z^{n_u-1} \\
%    \gamma_{i,q}(z)&=& \frac{1}{h_{q}^{n_u+n_y}}
%    (z^{n_u-i}-\sum_{j=1}^{i-1} \gamma_{j,q}(z)h_{q}^{n_y+n_u-i+j})
%  \end{eqnarray*}
  \end{Lemma}
  \begin{proof}%[Proof of Lemma \ref{sarx2lss:lemma3:lemma1}]
   In this proof we will view $A_{q}$ as a linear map
   $x^T \mapsto x^TA_{q}$, defined
  on the space of row vectors 
  $x^T \in \mathbb{R}^{1 \times (n_y+n_u)}$.
  From the structure of $A_{q}$ it then follows that
  $e_{j}^TA_{q}=e_{j-1}^T$, for $j=\{2,\ldots,n_y\}\cup \{n_y+2,\ldots,n_y+n_u\}$.
  Hence,  $e_{n_y}^TA_{q}^{j}$, $j=0,\ldots,n_y-1$
  spans $X_1=\mathrm{Span}\{e_1^T,\ldots,e_{n_y}^T\}$.
  Notice that
   $e_1^TA_{q}=\sum_{i=1}^{n_y+n_u-1} h_q^ie_i^T+h_q^{n_u+n_y}e_{n_y+n_u}^T$.
   Hence 
   \begin{equation}
  \label{sarx2lss:lemma3:lemma1:eq0}
   x^T=\sum_{i=n_y+1}^{n_y+n_u} h_{q}^{i}e^T_i=e_1^TA_{q}-\sum_{i=1}^{n_y} h_q^ie_i^T
  \end{equation}
  belongs to the
  linear span of $e_{n_y}^TA_{q}^{j}$, $j=0,\ldots,n_y$.

  We proceed to prove that $x^TA_{q}^{j}$, $j=0,\ldots,n_u-1$
  span
  $X_2=\mathrm{Span}\{e_{n_y+1}^T,\ldots e_{n_y+n_u}^T\}$.
  From this the first statement of the lemma follows.
  Notice that $e_{n_y+1}^TA_{q}=0$ and 
  $e_{n_y+j}^TA_q=e_{n_y+j-1}^T$ for all $j=2,\ldots,n_u$.
  Hence, 
  \begin{equation} 
  \label{sarx2lss:lemma3:lemma1:eq1}
    x^TA_{q}^{j}=\sum_{i=1}^{n_u-j} h_q^{n_y+i+j}e^T_{n_y+i}.
  \end{equation}
  From \eqref{sarx2lss:lemma3:lemma1:eq1} and
  $h_q^{n_y+n_u} \ne 0$ it then follows that
  \begin{equation}
  \label{sarx2lss:lemma3:lemma1:eq2}
    e_{n_y+1}^T=\frac{1}{h_q^{n_u+n_y}} x^TA_{q}^{n_u-1}
  \end{equation}
  and if $e_{n_y+1}^T,\ldots, e_{n_y+j}^T$ have already been
  obtained from the linear combinations of
  $x^TA_{q}^i$, $i = n_u-j, \ldots,n_u-1$, then
  \begin{equation}  
  \label{sarx2lss:lemma3:lemma1:eq3}
    e_{n_y+j+1}^T = \frac{1}{h_{q}^{n_u+n_y}} (x^TA_{q}^{n_u-j-1} - \sum_{i=1}^{j} h_{q}^{n_y+n_u-j-1+i}e_{n_y+i}).
  \end{equation}
  Hence, $x^TA_{q}^j$, $j=0,\ldots, n_u-1$ span 
  $X_2$.

  Finally, the statement $e_{j}^T=e_{n_y}^TA_{q}^{n_y-j}$,
  $j=1,\ldots,n_y$ follows from the definition of 
  $A_q$.  The statement that $e_{n_y+j}^T=e_{n_y}^T\chi(A_q)\gamma_{j,q}(A_q)$, $j=1,\ldots,n_u$ can be shown as follows.
   From \eqref{sarx2lss:lemma3:lemma1:eq0} it follows that
   $x^T=e_{n_y}^T\chi_q(A_q)$. From 
   \eqref{sarx2lss:lemma3:lemma1:eq2} and 
   \eqref{sarx2lss:lemma3:lemma1:eq3} it follows
   that $e_{n_y+j}^T=x^T\gamma_{j,q}(A_q)$ for all $j=1,\ldots,n_u$.
   Combining the above statements implies the second statement
   of the lemma.
  \end{proof}
  \begin{Lemma}
  \label{sarx2lss:lemma3:lemma2}
   Assume that $h_{q}^{n_u+n_y} \ne 0$.
   The characteristic polynomial of $A_{q}$ coincides with
   its minimal polynomial and it equals
   \[ \chi_{A_q}(z)=z^{n_u}(z^{n_y}-\sum_{i=1}^{n_y} h_{q}^iz^{n_y-i}).
   \]
  \end{Lemma}
 \begin{proof}%[Proof of Lemma \ref{sarx2lss:lemma3:lemma2}]
  From Lemma \ref{sarx2lss:lemma3:lemma1} it follows that
  $\mathbb{R}^{1 \times (n_u+n_y)}$ is a cyclic subspace
  with respect to the linear operator $\hat{A}_q:x^T \mapsto x^TA_{q}$
  \footnote{The definition of cyclic subspaces can be found in Section 4-Chapter VII of the book of Gantmacher~\cite{GantmacherBook}}.
  By Theorem 4-Chapter VII of the book of Gantmacher~\cite{GantmacherBook}, it then follows that the minimal polynomial of the
  linear operator $\hat{A}_q$ equals its characteristic polynomial and it is of degree $n_y+n_u$.
  Note that in the standard basis $e_1^T,\ldots,e_{n_u+n_y}^T$,
  the basis of the linear operator $\hat{A}_q$ is
  $A_{q}^T$. Hence, the minimal polynomial and characteristic
  polynomial of $A_{q}^T$ coincide. But these polynomials 
  are the same for the matrices $A_q$ and $A_q^T$.

  Moreover, from  Lemma \ref{sarx2lss:lemma3:lemma1}, it
  also follows that $e_{n_y}^T$ is the generating element
  of the cyclic space $\mathbb{R}^{1 \times (n_y+n_u)}$.
  Hence, by Subsection 4.1-Chapter VII of the book of Gantmacher~\cite{GantmacherBook},
  a polynomial $\psi(z)$ is a minimal polynomial of  $\hat{A}_q$, if $\psi(\hat{A}_q)e_{n_y}^T=e_{n_y}^T\psi(A_q)=0$ and
  it has the smallest possible degree. By the discussion above,
  the degree of the minimal polynomial of $\hat{A}_q$ must
  be $n_y+n_u$. Hence, the minimal polynomial of $\hat{A}_q$ is
  the unique monic polynomial $\psi(z)$ of degree $n_y+n_u$,
  such that $e_{n_y}^T\psi(A_q)=0$.

  %Consider now $\=z^{n_y+n_u}-\sum_{i=1}^{n_y} h_{q}^iz^{n_y-i+n_u})$.  
  If we show that $e_{n_y}^T\chi_{A_q}(A_q)=0$, then
  the statement of the lemma follows.
  To this end, notice that if 
  $y^T \in X_2=\mathrm{Span}\{e_{n_y+1}^T,\ldots,e^T_{n_y+n_u}\}$,
  then $y^TA^{n_u}_{q}=0$. In addition,
  $e_{n_y}^TA_{q}^{n_y}=e_1^TA_{q}=\sum_{i=1}^{n_y} h_{q}^{i}e_i^T+x^T$, where $x^T=\sum_{i=n_y+1}^{n_u} h_q^{i}e_{i}^T \in X_2$.
  Hence, by taking into account the remark above and 
  that $e_i^T=e_{n_y}A_{q}^{n_y-i}$, $i=1,\ldots,n_y$, 
  \begin{displaymath}
  e_{n_y}^TA_{q}^{n_y+n_u}=
  \sum_{i=1}^{n_y} h_{q}^{i}e_i^TA_{q}^{n_u}+x^TA_{q}^{n_u}=
  \sum_{i=1}^{n_y} h_{q}^{i}e_{n_y}^TA_{q}^{n_u+n_y-i}.
  \end{displaymath}
%  \[ 
%    \begin{split}
%     & e_{n_y}^TA_{q}^{n_y+n_u}=
%     \sum_{i=1}^{n_y} h_{q}^{i}e_i^TA_{q}^{n_u}+x^TA_{q}^{n_u}= \\
%     & \sum_{i=1}^{n_y} h_{q}^{i}e_{n_y}^TA_{q}^{n_u+n_y-i}
%    \end{split}.
%  \]
  The latter is exactly equivalent to
  $e_{n_y}^T\chi_{A_q}(A_q)=0$.  
 \end{proof}
%%%%%%%%%%%%%%%%%%%%%%%%%%%

Now, to complete the proof of Theorem \ref{sarx2lss:lemma3}, we will show that if Part \textbf{(A)} holds, then
 $\Sigma_{\AS}$ is reachable, and if Part \textbf{(B)}
 holds, then $\Sigma_{\AS}$ is observable.
% \begin{proof}[Proof Theorem \ref{sarx2lss:lemma3}]

 \textbf{Proof of Part (A)} \\
   We will show that if the conditions of \textbf{(A)} hold, then
   $(A_{q_0},A_{q_1}^{n_u}B_{q_1})$ is a controllable pair.
   By Sun and Ge~\cite{Sun:Book} it then follows that the \LSS\ 
   $\Sigma_{\AS}$ is reachable.
   From Lemma \ref{sarx2lss:lemma2:lemma2.1} it follows that 
   $A_{q_1}^{j}B_{q_1}=A_{q_1}^{j}e_{n_y+1}=h_{q_1}^{n_y+j}A_{q_1}^{j-1}e_1+e_{n_y+j+1}$ for $j=1,\ldots,n_u-1$ and hence
  \[ A_{q_1}^{n_u}B_{q_1}=\sum_{j=1}^{n_u} h_{q_1}^{n_y+j}A_{q_1}^{n_u-j}e_1. \]
  From Lemma \ref{sarx2lss:lemma2:lemma2.1} it also follows that
  $\psi_{q_0,q_1,j}(A_{q_0})=A_{q_1}^{j}$ and hence 
  the polynomial $\phi_{q_0,q_1}(z)$ satisfies
  \[ A_{q_1}^{n_u}B_{q}=\phi_{q_0,q_1}(A_{q_0})e_1.    \]
  From Lemma \ref{sarx2lss:lemma2:lemma2.1}, it follows that
  $\phi_{q_0,q_1}(A_{q_0})e_1 \in X_1$ and $X_1$ is $A_{q_1}$ invariant,
  where $X=\mathrm{Span}\{e_1,\ldots,e_{n_y}\}$. In addition,
  from the construction of $A_{q_0}$ it follows that with respect to
  the basis $e_1,\ldots,e_{n_y}$, the matrix representation of 
  the restriction of $A_{q_0}$ to $X_1$ is of the form
  \[
    \hat{A}_{q_0}=\begin{bmatrix}
     \begin{bmatrix} h_{q_0}^{1} & \ldots & h_{q_0}^{n_y-1} \end{bmatrix} & h_{q_0}^{n_y} \\
     I_{n_y-1} &  \mathbf{O}_{(n_y-1) \times 1} 
    \end{bmatrix}.
  \]
  The above matrix is in companion form and it is known that its
  characteristic polynomial equals its minimal polynomial and it
  equals $\chi_{q_0}(z)$. That is, $\chi_{q_0}(z)$ is the minimal
  polynomial of the linear operator $A_{q_0}$ restricted to $X_1$.
  Moreover, from  Lemma \ref{sarx2lss:lemma2:lemma2.1}, it follows
  that $A_{q_0}^{j}e_1$, $j=0,\ldots, n_{y}-1$ generate the space
  $X_1$, i.e. $X_1$ is a cyclic subspace w.r.t. to $A_q$.
  Then by Subsection 4.1--Chapter VII of the book of Gantmacher~\cite{GantmacherBook}, 
  $\chi_{q_0}(z)$ is a minimal polynomial of $e_1$ with respect
  to $A_{q_0}$, i.e. $\chi_{q_0}(A_{q_0})e_1=0$ and 
  $\chi_{q_0}(z)$ has the smallest degree among all the polynomials
  $\psi(z)$ such that $\psi(A_{q_0})e_1=0$.
  
  Suppose now that $\chi_{q_0}(z)$ and $\phi_{q_0,q_1}(z)$ are coprime, but
  $(A_{q_0},A_{q_1}^{n_u}B_{q_1})=(A_{q_0},\phi_q(A_{q_0})e_1)$ is not a
  controllable pair.  Then the vectors 
  $A_{q_0}^{j}x$, $j=0,\ldots,n_y-1$,
  $x=\phi_{q_0,q_1}(A_{q_0})e_1$ are linearly dependent, i.e. there exists
  a non-zero polynomial $\kappa(z)$ of degree at most $n_y-1$ such
  that $\kappa(A_{q_0})x=0$. By substituting $x=\phi_{q_0,q_1}(A_{q_0})e_1$,
  we get $\kappa(A_{q_0})\phi_{q_0,q_1}(A_{q_0})e_1=0$. That is, for the
  polynomial $\phi(z)=\kappa(z)\phi_{q_0,q_1}(z)$, 
  $\phi(A_{q_0})e_1=0$. This implies by Gantmacher~\cite{GantmacherBook} that the minimal polynomial $\chi_{q_0}(z)$ 
  divides $\phi(z)=\kappa(z)\phi_{q_0,q_1}(z)$. Since $\chi_{q_0}(z)$ and
  $\phi_{q_0,q_1}(z)$ are co-prime, then this is possible only if 
  $\chi_{q_0}(z)$ divides $\kappa(z)$. But the degree of $\kappa(z)$ is
  strictly smaller than the degree of $\chi_{q_0}(z)$, hence 
  $\kappa(z)$ cannot be divisible by $\chi_{q_0}(z)$. We arrived to a 
  contradiction. That is, we can conclude that
  $(A_{q_0},A_{q_1}^{n_u}B_{q_1})$ is a controllable pair.

  \textbf{Proof of Part \textbf(B)} \\
  We will show that $(C_{q_3},A_{q_2})$ is an observable pair.
  By Sun and Ge~\cite{Sun:Book}, this is sufficient for
  observability of $\Sigma_{\AS}$.  

  To this end, using the notation of Lemma \ref{sarx2lss:lemma3:lemma1}
  define the polynomial
  \[ \hat{\psi}(z)=\upsilon_{q_3}(z)+\sum_{j=1}^{n_u} h_{q_3}^{n_y+j}
      \gamma_{j,q_2}(z)\chi_{q_2}(z).
  \]
   Then from Lemma \ref{sarx2lss:lemma3:lemma1} it follows that $C_{q_3}=e_{n_y}^T\hat{\psi}(A_{q_2})$.
  Assume that $(C_{q_3},A_{q_2})$ is not an observable pair.
  Then $C_{q_3}A_{q_2}^{j}$, $j=0,\ldots,n_y-1$ are
  linearly independent. Hence, there exists a polynomial
  $\kappa(z)$ of degree less than $n_y$, such that 
  $C_{q_3}\kappa(A_{q_2})=0$. Hence, we obtain that
  $e_{n_y}^T\hat{\psi}(A_{q_2})\kappa(A_{q_2})=0$.
  In other words, the polynomial $P(z)=\hat{\psi}(z)\kappa(z)$
  is an annihilating polynomial with respect to the operator
  $\hat{A}_{q_2}:x \mapsto xA_{q_2}$ \footnote{For the definition, see the book of Gantmacher~\cite{GantmacherBook}.} of $e_{n_y}^T$.  
Since by Lemma \ref{sarx2lss:lemma3:lemma1} $e_{n_y}^TA_{q_2}^{j}$, $j=0,\ldots,n_y+n_u$ generate
  the whole space, it then follows that $P(z)$ is the annihilating polynomial of the whole space, i.e. $P(A_{q_2})=0$. It then
  follows that $P(z)$ is divisible by the minimal polynomial
  of $\hat{A}_{q_2}$ which coincides with that of $A_{q_2}$.
  From Lemma \ref{sarx2lss:lemma3:lemma2} it follows that
  the minimal polynomial of $A_{q_2}$ is $z^{n_u}\chi_{q_2}(z)$.

  We will argue that if the conditions of Part \textbf{(B)}
  hold, then $\hat{\psi}(z)$ and $z^{n_u}\chi_{q_2}(z)$ are
  co-prime. Indeed, if $\hat{\psi}(z)$ and $z^{n_u}\chi_{q_2}(z)$
  are not co-prime, then there exists an irreducible polynomial
  $q(z)$ which divides both $\hat{\psi}(z)$ and $z^{n_u}\chi_{q_2}(z)$.
  If $q(z)$ is an irreducible polynomial which divides
  $z^{n_u}\chi_{q_2}(z)$, then it either equals $z$ or it divides
  $\chi_{q_2}(z)$. If $q(z)=z$ and it divides $\hat{\psi}(z)$, then
  $0$ is a root of $\hat{\psi}(z)$, i.e. 
  $\hat{\psi}(0)=0$. Notice that by induction it follows that
  for $j=1,\ldots,n_u-1$, $\gamma_{j,q_2}(0)=0$ and 
  $\gamma_{n_u,q_2}(0)=\frac{1}{h_{q_2}^{n_u+n_y}}$. Hence, from
  the definition of $\hat{\psi}(z)$ it follows that
  $\hat{\psi}(0)=h_{q_3}^{n_y}+\frac{h_{q_3}^{n_u+n_y}}{h_{q_2}^{n_u+n_y}}\chi_{q_2}(0)=h_{q_3}^{n_y}-\frac{h_{q_3}^{n_u+n_y}}{h_{q_2}^{n_u+n_y}}h_{q_2}^{n_y}$.
  Hence, $\hat{\psi}(0)=0$ implies that
  $h_{q_3}^{n_y}=\frac{h_{q_3}^{n_u+n_y}}{h_{q_2}^{n_u+n_y}}h_{q_2}^{n_y}$, which contradicts to the condition of \textbf{(B)}.
  If $q(z)$ divides $\chi_{q_2}(z)$ and it divides $\hat{\psi}(z)$,
  then it divides
  $\upsilon_{q_3}(z)=\hat{\psi}(z)-(\sum_{i=1}^{n_u} h_{q_3}^{n_u+i}\gamma_{i,q_2}(z))\chi_{q_2}(z)$. But this contradicts to the assumption
  that $\psi_{q_3}(z)$ and $\chi_{q_2}(z)$ are co-prime.

  Hence, by the discussion above, $\hat{\psi}(z)$ and
  $z^{n_u}\chi_{q_2}(z)$ are coprime, so if $z^{n_u}\chi_{q_2}(z)$
  divides $P(z)$, it then must divide $\kappa(z)$. But the
  degree of $\kappa(z)$ is strictly smaller than that of
  $z^{n_u}\chi_{q_2}(z)$, hence $z^{n_u}\chi_{q_2}(z)$ cannot divide 
  $\kappa(z)$. We arrived to a contradiction. Hence,
  $(C_{q_3},A_{q_2})$ must be an observable pair.

 %\end{proof}
%%%%%%%%%%%%%%
%\section{Poof of Corollary~\ref{sarx2lss:col2}}\label{appdxProofCorosarx}
%\section{Proof of Theorem \ref{sarx:theo2}}\label{appxTheo4}
\section{Proof of Theorem \ref{sarx:main2.1}}\label{appxTheo3}
In order to prove Theorem \ref{sarx:theo2} we will relate identifiability analysis of SARX systems to that of the associated \SLSS\,\ntext{(see Section~\ref{sect:ident} for the definition of a parametrization and identifiability of \SLSS)}. This is possible due to the following corollary of Lemma \ref{sarx2lss:lemma0}.
\begin{Corollary}
\label{sarx2lss:col2}
 A SARX parametrization $\MAPSARX$ of the form \eqref{sarx:param:def:eq1} is identifiable, if and only if
 the \LSS\ parametrization $\mathbf{\Pi}_{sw}: \PARAMSET \ni \theta \mapsto \Sigma_{\MAPSARX(\theta)}$ is identifiable. \ntext{Here, $\Sigma_{\MAPSARX(\theta)}$ is the LSS of the form \eqref{sarx2lss:eq1}--\eqref{sarx2lss:eq1bis} obtained from the SARX $\MAPSARX(\theta)$.}
\end{Corollary}
\begin{proof}[Proof of Corollary~\ref{sarx2lss:col2}]
\color{black}
Consider two SARX systems $\AS_i=\{h_q^i\}_{q \in Q}$,
  $i=1,2$ of type $(n_y,n_u)$.
  Notice that each $\AS_i$, $i=1,2$,  realizes the same input-output map as the associated \SLSS\ $\Sigma_{\AS_i}$.\\ 
  %$\Sigma_{\AS_2}$ realize the same input-output map. \\
%%%%%
  Assume that the parametrization $\MAPSARX$ is
  identifiable, but $\mathbf{\Pi}_{sw}$ is not identifiable.
  Then there exist two parameters $\theta_1,\theta_2 \in \PARAMSET$, $\theta_1 \ne \theta_2$, such that
  $\mathbf{\Pi}_{sw}(\theta_1)$ and $\mathbf{\Pi}_{sw}(\theta_2)$
  realize the same input-output map. Since
  $\mathbf{\Pi}_{sw}(\theta_i)=\Sigma_{\MAPSARX(\theta_i)}$,
  $i=1,2$, by the remark above it follows that
  $\MAPSARX(\theta_1)$ and $\MAPSARX(\theta_2)$ are
  equivalent. This contradicts the identifiability of
  $\MAPSARX$.

  Conversely, assume that $\mathbf{\Pi}_{sw}$ is identifiable,
  but $\MAPSARX$ is not identifiable. Then 
  there exists parameters $\theta_1,\theta_2 \in \PARAMSET$,
  $\theta_1 \ne \theta_2$, such that
  $\MAPSARX(\theta_1)$ and $\MAPSARX(\theta_2)$ are
  equivalent. This means that 
  $\Sigma_{\MAPSARX(\theta_1)}=\mathbf{\Pi}_{sw}(\theta_1)$ and
  $\Sigma_{\MAPSARX(\theta_2)}=\mathbf{\Pi}_{sw}(\theta_2)$ 
  realize the same input-output map. But this contradicts the
  identifiability of $\mathbf{\Pi}_{sw}$.
\end{proof}
  In order to prove Theorem \ref{sarx:main2.1}, we will
  need the following result which is interesting on its
  own right.
\begin{Theorem}
 \label{sarx:theo2}
  Consider two SISO SARX systems $\AS_1=\{h_q\}_{q \in Q}$ and
  $\AS_2=\{g_q\}_{q \in Q}$ of type $(n_y,n_u)$ and assume that for some $q \in Q$, either $h_q^{n_y} \ne 0$ or $h_{q}^{n_y+n_u} \ne 0$.
  %Then the only isomorphism which exists between
  %the associated \SLSS\ $\Sigma_{\AS_1}$ and 
  %$\Sigma_{\AS_2}$ is the identity map.
   If there exists an isomorphism\footnote{See the work of Petreczky {\it et al}~\cite{MP:HSCC2010} for the definition of isomorphism between \SLSS} 
  between
  the associated \SLSS\ $\Sigma_{\AS_1}$ and 
  $\Sigma_{\AS_2}$, then this isomorphism is the identity map.
 \end{Theorem} 
  Theorem \ref{sarx:theo2} implies,  under some
  mild conditions, that the transformation
  of two different SARX systems 
  to state-space representations cannot result in 
  isomorphic systems.
 \begin{proof}[Proof of Theorem \ref{sarx:theo2}]
\color{black}
  Assume that $\Sigma_{\AS_1}=\SwitchSysLin$ and $\Sigma_{\AS_2}=\SwitchSysLin[']$ with $n=n^{'}=n=n_y+n_u$ and $p=m=1$.
  Consider an isomorphism $\MORPH$ between $\Sigma_{\AS_1}$ and  
  $\Sigma_{\AS_2}$. Denote by $e_i$ the $i$th standard unit vector
  of $\mathbb{R}^{n}$. 
   Then $e_1^T,\ldots,e_n^T$
   form the standard  basis in $\mathbb{R}^{1 \times n}$. 
  The  proof depends on the following series of
  technical results.
  \begin{Proposition}
  \label{sarx:theo2:pf:prop0}
  \begin{equation}
  \label{sarx:theo2:pf:eq1}
     e_1^T A_{q} = e_1^T \MORPH A_{q}
  \end{equation}
  \end{Proposition}
  \begin{proof}%[Proof of Proposition \ref{sarx:theo2:pf:prop0}]
   From the construction of $\Sigma_{\AS_i}$, $i=1,2$ it then follows
   that $C_{q}=e_1^TA_{q}$, $C_{q}^{'}=e_1^TA_{q}^{'}$. From the definition of
   isomorphism between \SLSS, it follows that
   $C^{'}_{q}\MORPH=C_{q}$, $q \in Q$. Hence, we obtain that
   \[ e_1^TA_{q}=e_1^TA_{q}^{'}\MORPH. \]
   But $A_{q}^{'} \MORPH = \MORPH A_{q}$ by the definition of
   a \LSS\ isomorphism, and hence we obtain the claim of the
  proposition.
  \end{proof}
  \begin{Proposition} 
  \label{sarx:theo2:pf:prop1}
  The columns of $A_q$ span the space
  \[ \mathrm{Span}\{e_1,\ldots,e_{n_u+n_y}\}\setminus \{e_{n_y+1}\}. \]
  \end{Proposition}
  \begin{proof}[Proof of Proposition \ref{sarx:theo2:pf:prop1}]
   Indeed, 
  $A_{q}e_{n}=h_{q}^{n}e_1$, $A_{q}e_{n_y}=h_{q}^{n_y}e_1$,
  $A_{q}e_j=e_{j+1}+h_{q}^{j}e_1$
  for all $j \in \{1,\ldots,n-1\} \setminus \{n_y\}$. Hence, 
  if either $h_{q}^{n_y} \ne 0$ or $h_{q}^{n} \ne 0$, then
  $e_1$ belongs to the column space of $A_q$, and hence
  $e_j=A_{q}e_{j-1}-h_{q}^{j-1}e_1$ belongs to the column space
  of $A_q$, for $j\in \{2,\ldots,n\} \setminus \{n_y+1\}$.
  \end{proof}
 \begin{Proposition}
 \label{sarx:theo2:pf:prop2}
  For any $i=1,\ldots,n_u+n_y$,
  if $e_{i}^TA_{q}=e_i^T\MORPH A_{q}$, then 
  $e_i^T=e_i^T\MORPH$.
 \end{Proposition}
 \begin{proof}[Proof of Proposition \ref{sarx:theo2:pf:prop2}]
  Indeed, if $e_i^TA_{q}=e_i^T\MORPH A_{q}$, then
  this implies that $(e_i^T-e_i^T\MORPH)A_{q}=0$.
  By Proposition \ref{sarx:theo2:pf:prop1} this implies that
  $(e_i^T-e_i^T\MORPH)e_j=0$ for all $j \in \{1,\ldots,n_u+n_y\}\setminus \{n_y+1\}$.
  Notice that  from the construction of
  $\Sigma_{\AS_1}$, $\Sigma_{\AS_2}$ and the
  definition of a \LSS\ morphism it follows that
  $e_{n_y+1}=B_q^{'}=\MORPH B_{q} =\MORPH e_{n_y+1}$. Hence,
  $(e_i^T-e_i^T\MORPH)e_{n_y+1}=0$ and thus
  \[ (e_i^T-e_i^T\MORPH)e_j=0, j=1,\ldots,n_y+n_u. \]
 This is just an alternative way of formulating the
 conclusion of the proposition.
 \end{proof}
%%%
 \begin{Proposition}
 \label{sarx:theo2:pf:prop3}
  If $e_{j-1}^T\MORPH=e_{j-1}^T$, then
  $e_{j}^TA_{q}=e_j^T\MORPH A_q$ for all
  $j=\{2,\ldots,n_y+n_u\} \setminus \{n_y+1\}$.
 \end{Proposition}
 \begin{proof}%[Proof of Proposition \ref{sarx:theo2:pf:prop3}]
  Notice that $e_j^TA_{q}=e_{j-1}^T$, 
  $e_j^TA_{q}^{'}=e_{j-1}^T$,
  $j=n_y,\ldots,2$, and
  $e_{j}^TA_{q}=e_{j-1}^T$, $e_j^TA_{q}^{'}=e_{j-1}^T$, 
  for $j=n_{y}+n_{u},\ldots, n_{y}+2$.  
   Hence, by using $A_{q}^{'}\MORPH=\MORPH A_{q}$, we derive
  \begin{equation}
  \label{sarx:theo2:pf:eq2}
     e_{j-1}^T\MORPH = e_{j}^TA_{q}^{'}\MORPH=e_{j}^T\MORPH A_{q}
  \end{equation}
  for all $j \in \{2,\ldots,n_y\} \cup \{n_y+2,\ldots,n_{y}+n_u\}$. Since $e_{j-1}^T\MORPH = e_{j-1}^T$, and
  $e_j^TA_{q}=e_{j-1}^T$ for all
  $j=\{2,\ldots,n_u+n_y\}\setminus\{n_y+1\}$, from
  \eqref{sarx:theo2:pf:eq2} we obtain
  the claim of the proposition.
 \end{proof}
The rest of the proof of Theorem \ref{sarx:theo2} proceeds as follows.
 We will prove that 
 \begin{equation}
 \label{sarx:theo2:pf:eq3}
  e_j^T=\MORPH e_j^T, j=1,\ldots,n_y+n_u,   
 \end{equation}
 which is just another way of saying that $\MORPH$ is the
 identity matrix.
 To this end, from
 \eqref{sarx:theo2:pf:eq2} and Proposition
 \ref{sarx:theo2:pf:prop2} it follows that
 \eqref{sarx:theo2:pf:eq3} holds for $j=1$.
 Moreover, the $n_y+1$th row of $A_{q}$ and $A_{q}^{'}$ are
 both zero, hence, 
 $0=e_{n_y+1}^TA_{q}$ $0=e_{n_y+1}^TA_{q}^{'}$ and
 thus $0=e_{n_y+1}A_{q}^{'}\MORPH=e_{n_y+1}^T\MORPH A_{q}$.
 From this we get that $e_{n_y+1}^TA_{q}=e_{n_y+1}^T\MORPH A_{q}$ and by Proposition \ref{sarx:theo2:pf:prop2} this 
 implies that \eqref{sarx:theo2:pf:eq3} holds for $j=n_y+1$.
 Notice that if \eqref{sarx:theo2:pf:eq3} holds for
 $j=k \in \{1,\ldots,n_u+n_y-1\}\setminus \{n_y\}$, then
 by Proposition \ref{sarx:theo2:pf:prop3}, 
 $e_{k+1}^TA_{q}=e_{k+1}^T\MORPH A_{q}$. By Proposition 
 \ref{sarx:theo2:pf:prop2}, the latter implies 
 that \eqref{sarx:theo2:pf:eq3} holds for $j=k+1$.
 Hence, by induction we get that
 \eqref{sarx:theo2:pf:eq3} holds for all $j$.
\end{proof}
%%%%%%%%%%%%%%%%%%%%%%%
 \begin{proof}[Proof of Theorem \ref{sarx:main2.1}]
 We will show that 
  the \LSS\ parameterization $\mathbf{\Pi}_{sw}:\PARAMSET \ni \theta \mapsto \Sigma_{\MAPSARX(\theta)}$  is identifiable.
  By Corollary \ref{sarx2lss:col2} this is sufficient
  for identifiability of $\MAPSARX$.

  Since $\MAPSARX$ is strongly minimal, 
  the \LSS\ parameterization $\mathbf{\Pi}_{sw}$ is
  minimal~\cite{MP:HSCC2010}.
  In order to show identifiability of $\mathbf{\Pi}_{sw}$, 
  according to Petreczky {\it et al}~\cite[Corollary 1]{MP:HSCC2010}, it is enough to show
  that the only isomorphism between elements of 
  $\mathbf{\Pi}_{sw}$ is the identity.
  Consider now two elements
  $\Sigma_i=\Sigma_{\MAPSARX(\theta_i)}$, $\theta_i \in \PARAMSET$, $i=1,2$ of $\mathbf{\Pi}_{sw}$.
  Notice that $\MAPSARX(\theta_1)$ is minimal, since it is strongly minimal, and  thus if $\MAPSARX(\theta_1)=\{h_q\}_{q \in Q}$, then by Lemma \ref{element:lemma1} $h_{q}^{n_u+n_y} \ne 0$.
  But then from Theorem \ref{sarx:theo2} it follows that the only isomorphism between $\Sigma_1$ and $\Sigma_2$ is the identity map.

%%  From strong minimality of $\MAPSARX(\theta)$
%%  it follows that $\AS=\MAPSARX(\theta)$ is
%%  minimal and hence if $\AS=\{h_q\}_{q \in Q}$, then
%%  by Lemma \ref{element:lemma1} 
%%   for some $q \in Q$, $h_q^{n_u+n_y} \ne 0$. 
%%  Hence, by Theorem \ref{sarx:theo2},
%%  for any $\theta_1 \ne \theta_2 \in G$, the only 
%%  \LSS\ isomorphism between 
%%  $\Sigma_{\MAPSARX(\theta_1)}$ and 
%%  $\Sigma_{\MAPSARX(\theta_2)}$ is the identity. But the 
%%   latter implies that $\Sigma_{\MAPSARX(\theta_1)}=\Sigma_{\MAPSARX(\theta_2)}$ which implies that $\MAPSARX(\theta_1)=\MAPSARX(\theta_2)$. Since $\MAPSARX$ is an injective map, 
%%  $\theta_1=\theta_2$ which is a 
%%  contradiction. That is, there exists no isomorphism between
%%  $\Sigma_{\MAPSARX(\theta_1)}$ and $\Sigma_{\MAPSARX(\theta_2)}$.
%%  Hence, by \cite{MP:HSCC2010} the \LSS\ parametrization
%%  $\mathbf{\Pi}_{sw}:\PARAMSET \ni \theta \mapsto \Sigma_{\MAPSARX(\theta)}$ is identifiable.
\end{proof}
\section{Proof of Theorem \ref{sarx:theo1}}\label{appxTheo5}
 %\begin{proof}[Proof of Theorem \ref{sarx:theo1}]
  Let $K=(pn_y+mn_u)|Q|$. Then any SARX system of type $(n_y,n_u)$
  can be identified with a point in $\mathbb{R}^{K}$, by identifying
  the system with the collection of its parameters $\{h_q\}_{q \in Q}$, $h_q \in \mathbb{R}^{p \times (pn_y+mn_u)}$. 
  
  First, we construct a polynomial 
  $P_{min}(X_1,\ldots,X_{K})$,
  such that $P_{min}(\AS) \ne 0$ if and only if $\AS$ is strongly
  minimal. To this end, consider the \LSS\ $\Sigma_{\AS}$ and
  consider the observability and controllability matrices 
  $O(\Sigma_{\AS})$ and $\mathcal{R}(\Sigma_{\AS})$ as defined in \cite{Petreczky2012-Automatica}. 
  %%%%%%%%%%%%%%%%%
  Define 
  \begin{displaymath}
  P_{obs}(X_1,\ldots,X_K)=\det(O(\Sigma_{\AS})^T O(\Sigma_{\AS})),
  \end{displaymath}
  and
  \begin{displaymath}
  P_{contr}(X_1,\ldots,X_K)=\det(\mathcal{R}(\Sigma_{\AS})\mathcal{R}(\Sigma_{\AS})^T),
  \end{displaymath}
%  $$ \begin{aligned}
%    &P_{obs}(X_1,\ldots,X_K)=\det(O(\Sigma_{\AS})^T O(\Sigma_{\AS}))\\
%     &P_{contr}(X_1,\ldots,X_K)=\det(\mathcal{R}(\Sigma_{\AS})\mathcal{R}(\Sigma_{\AS})^T),
%\end{aligned}$$ 
with the notation $\det(\cdot)$ referring to determinant. Then $P_{obs}$ and $P_{contr}$ are polynomials in the entries of the matrices of $\Sigma_{\AS}$.
Moreover, by applying the results of Petreczky {\it et al}~\cite{Petreczky2012-Automatica}, 
\begin{itemize}
	\item $P_{obs}(\AS) \ne 0$ $\Leftrightarrow$ $\mathcal{O}(\Sigma_{\AS})$ has full rank $\Leftrightarrow$ $\Sigma_{\AS}$ is observable
\item $P_{contr}(\AS) \ne 0$ $\Leftrightarrow$   $\mathcal{R}(\Sigma_{\AS})$ has full rank $\Leftrightarrow$ $\Sigma_{\AS}$ is reachable.
\end{itemize}
  Define now \( P_{min}=P_{contr}P_{obs}  \).
  Then $P_{min}(\AS) \ne 0$ if and only if $\Sigma_{\AS}$ is both
  observable and reachable, i.e. if and only if $\Sigma_{\AS}$ is
  minimal. 
  
  %Any $n \times n$
  %minor of $O(\Sigma_{\AS})$ or of $\mathcal{R}(\Sigma_{\AS})$ can
  %be viewed as a polynomial in the entries of the matrices
  %of $\Sigma_{\AS}$. Since the entries of the matrices of 
  %$\Sigma_{\AS}$ are linear functions of the entries of
  %the parameters of $\AS$, it follows that any $n \times n$ minor
  %of $O(\Sigma_{\AS})$ and $\mathcal{R}(\Sigma_{\AS})$ can
  %be viewed as  polynomial in $\AS$, where $\AS$ is identified with
  %an element of $\mathbb{R}^{K}$. Let $\mathcal{P}_1$ and
  %$\mathcal{P}_{2}$ be the
  %set of all $n \times n$ minors of $O(\Sigma_{\AS})$ 
  %and respectively
  %$R(\Sigma_{\AS})$, each minor being
  %viewed as a polynomial  $\mathbb{R}^{K}$.
  %Define the polynomials
  %\[
    %\begin{split}
     %& P_{obs}(X_1,\ldots,X_K)=\sum_{P \in \mathcal{P}_1} (P(X_1,\ldots,X_K))^{2}  \\
     %& P_{contr}(X_1,\ldots,X_K)=\sum_{P \in \mathcal{P}_2} (P(X_1,\ldots,X_K))^2  \\
    %\end{split}
  %\] 
  %It is then clear that $P_{obs}(\AS) \ne 0$ if and only if
  %at least one of the $n \times n$ minors of $O(\Sigma_{\AS})$ is
  %not zero, i.e. if and only if $\Rank O(\Sigma_{\AS})=n$.
  %That is, $P_{obs}(\AS) \ne 0$ if and only if $\Sigma_{\AS}$
  %is observable.
  %Similarly, $P_{contr}(\AS) \ne 0$ if and only if
  %$\Rank \mathcal{R}(\Sigma_{\AS}) \ne 0$, i.e if and only if
  %$\Sigma_{\AS}$ is reachable.
  %Define now
  %\( P_{min}=P_{contr}P_{obs}  \).
  %Then $P_{min}(\AS) \ne 0$ if and only if $\Sigma_{\AS}$ is both
  %observable and reachable, i.e. if and only if $\Sigma_{\AS}$ is
  %minimal. 

  Finally, consider a polynomial parametrization $\MAPSARX$
  such that $\MAPSARX$ contains a strongly minimal element.
  The fact that $\MAPSARX$ is a polynomial parametrization
  implies that there exists polynomials
  $\mathbf{\Pi}_i^{\text{\tiny \textbf{SARX}}}$ in variables $X_1,\ldots,X_d$, $i=1,\ldots, K$
  such that 
  $\MAPSARX(\theta)=(\mathbf{\Pi}_1^{\text{\tiny \textbf{SARX}}}(\theta),\ldots,\mathbf{\Pi}_K^{\text{\tiny \textbf{SARX}}}(\theta))$
  for all $\theta \in \PARAMSET$. Here we used the identification of
  a SARX system of type $(n_y,n_u)$ with a point in
  $\mathbb{R}^{K}$.
  Consider the polynomial
  \[ 
  \begin{split}
    Q_{min}(X_1,\ldots,X_d)=P_{min}(\mathbf{\Pi}_1^{\text{\tiny \textbf{SARX}}}(X_1,\ldots,X_d),\ldots
   \ldots, \mathbf{\Pi}_K^{\text{\tiny \textbf{SARX}}}(X_1,\ldots,X_d)).
  \end{split}
  \]
  Notice that the set of parameters from $\PARAMSET$ which do not
  yield a minimal SARX system all satisfy the equation
  $Q_{min}(\theta)=0$. 
  From the assumption that $\MAPSARX$ contains a strongly minimal
  element it follows that for some $\theta \in \PARAMSET$,
  $Q_{min}(\theta)=P_{min}(\MAPSARX(\theta)) \ne 0$. Hence, the set
  $G=\{ \theta \in \PARAMSET \mid Q_{min}(\theta) \ne 0 \}$ is a 
  non-empty subset of $\PARAMSET$ and it is clearly generic.
  That is, $\MAPSARX$ is generically strongly minimal, and hence
  minimal.
% \end{proof}
%%%
 %\begin{proof}[Proof of Corollary \ref{sarx:col_main2}]
  %If $\MAPSARX$ is generically strongly minimal, then there
  %exists a generic set $G \subseteq \PARAMSET$ such
  %that the parametrization 
  %$\MAPSARX|_{G}:G \ni \theta \mapsto \MAPSARX(\theta)$
  %is strongly minimal. Hence, by Theorem \ref{sarx:main2.1},
  %$\MAPSARX|_{G}$ is identifiable. This means that 
  %$\MAPSARX$ is generically strongly identifiable.
 %\end{proof}
%%%%
 %\begin{proof}[Proof of Corollary \ref{sarx:main2}]
  %If $\MAPSARX$ contains a strongly minimal element,
  %then by Theorem \ref{sarx:theo1} $\MAPSARX$ is
  %generically strongly minimal. The rest follows from
  %Corollary \ref{sarx:col_main2}.
 %\end{proof}
%%%%%%%%%%%%%%%%%%%%%
 %\begin{proof}[Proof of Corollary \ref{sarx:main:col2}]
  %By Example \ref{example:min2}, there exists a strongly minimal
  %SARX system, i.e. $SARX_{triv}$ contains a strongly minimal
  %element. Moreover, $SARX_{triv}$ is clearly injective and polynomial. The statement follows now
   %Theorem \ref{sarx:theo1} and Corollary
 %\ref{sarx:main2}.
%\end{proof}
%%%%%%%%%%%%%%%%%%%%%%%%%%%%%%%%%%%%%%ù
%\subsection{Vehicle suspension system}\label{AppendxVS}

 %Note the converse is not true, since in general it is not
 %clear how to construct a SARX systems from a \LSS. We will
 %present a counter-example later.
%That is, minimality of the state-space representations of
%SARX systems cannot be deduced from linear theory.
%==============================================

%\nocite{*}% Show all bib entries - both cited and uncited; comment this line to view only cited bib entries;
%\bibliography{sysid2012_hybrid}%
\bibliography{Biblio-IJRNC-MP}
%\clearpage

%\section*{Author Biography}

%\begin{biography}{\includegraphics[width=66pt,height=86pt,draft]
%{empty}}{\textbf{Author Name.} This is sample author biography text this is sample author biography text this is sample author biography text this is sample author biography text this is sample author biography text this is sample author biography text this is sample author biography text this is sample author biography text this is sample author biography text this is sample author biography text this is sample author biography text this is sample author biography text this is sample author biography text this is sample author biography text this is sample author biography text this is sample author biography text this is sample author biography text this is sample author biography text this is sample author biography text this is sample author biography text this is sample author biography text.}
%\end{biography}

\end{document}